\documentclass[10pt,a4paper,twoside,notitlepage]{article}

\usepackage{graphicx}

\usepackage{amsfonts}
\usepackage{amsmath}
\usepackage{amssymb}
\usepackage{psfrag}
\usepackage{latexsym, epsfig}
\usepackage[active]{srcltx}
\usepackage{color}

\usepackage{wasysym}

\newtheorem{theorem}{Theorem}
\newtheorem{lemma}{Lemma}[section]
\newtheorem{corollary}{Corollary}[section]
\newtheorem{proposition}{Proposition}[section]
\newtheorem{definition}{Definition}[section]
\newtheorem{remark}{Remark}[section]

\newcommand{\R}{\mathbb{R}}
\newcommand{\ve}{\varepsilon}

\newcommand{\n}{\noindent}

\newcommand{\N}{\mathbb{N}}

\newcommand{\BV}{\mathbf{BV}}
\newcommand{\Lsp}{\mathbf{L}}

\newcommand{\diam}{\mathrm{diam}}

\newcommand{\A}{\mathcal{A}}

\newcommand{\ch}[1]{\mathrm{co}(#1)}
\newcommand{\chc}[1]{\overline{\mathrm{co}}\left(#1\right)}

\begin{document}

\title{
{\huge State constrained patchy feedback stabilization}
}

\author{
  \vspace{1cm}
    {\scshape Fabio S. Priuli}
    \thanks{Dipartimento di Matematica, Universit\`a degli Studi di Roma Tor Vergata, 
  Via della Ricerca Scientifica 1, I--00133 Roma, Italy.;                       
      e-mail: \texttt{priuli@mat.uniroma2.it}}     
      }

\date{\today}
\maketitle
\vspace{-0.10cm}
\begin{abstract}
\noindent
We construct a patchy feedback for a general control system on $\R^d$ which realizes practical stabilization to a target set $\Sigma$, when the dynamics is constrained to a given set of states $S$. The main result is that $S$--constrained asymptotically controllability to $\Sigma$ implies the existence of a discontinuous practically stabilizing feedback. Such a feedback can be constructed in ``patchy'' form, a particular class of piecewise constant controls which ensure the existence of local Carath\'eodory solutions to any Cauchy problem of the control system and which enjoy good robustness properties with respect to both measurement errors and external disturbances. 
\end{abstract}
\vspace{0.5cm}

2000\textit{\ Mathematical Subject Classification:} 34A; 49E; 93D


\vspace*{1cm}

\section{Introduction}
\label{sec:int}

Consider a general control system
\begin{equation}\label{eq:system}
\dot x = f(x,u)\qquad\quad x\in \R^d\,,~
\end{equation}
where  the upper dot denotes a derivative w.r.t.~time, $u$ is the control taking values in a compact set ${\bf U}\subset\subset\R^m$ and $f\colon\R^d\times{\bf U}\to\R^d$ is a vector field satisfying the following properties
\begin{description}
\item{{\bf (F1)}} $f$ is continuous on $\R^d\times{\bf U}$ and Lipschitz continuous in the variable $x$, uniformly for $u\in{\bf U}$, i.e. there exists a constant $L_f$ such that 
$$
|f(x,u)-f(y,u)|\leq L_f |x-y|\,,
$$
for all $(x,u)$ and $(y,u)$ in $\R^d\times{\bf U}$.
\item{{\bf (F2)}} $f$ has sub--linear growth, i.e.\begin{equation}\label{sublg}
|f(x,u)|\leq C_f \big(1+|x|\big)\qquad\quad\forall ~x\in \R^d\,,
\end{equation}
for some constant $C_f$ independent on $u$.
\item{{\bf (F3)}} The set of velocities
$$
f(x,{\bf U})\doteq\{f(x,u)~;~u\in{\bf U}\}
$$
is convex for every $x\in\R^d$.
\end{description}

\smallskip

In this article, we want to address the problem of stabilizing trajectories of~\eqref{eq:system} towards a target set $\Sigma\subseteq\R^d$ in the case of a dynamics constrained inside a prescribed set $S\subset \R^d$. In particular, we aim to construct a feedback control $U=U(x)$ which realizes stabilization to a neighborhood of $\Sigma$ and which is robust enough to provide the same stabilization also in presence of inner and outer perturbations of the dynamics, such as measurement errors and external disturbances.

\n However, one has to be careful because, even in very simple problems without state constraints, one cannot expect the existence of continuous control feedback laws which steer all trajectories towards a target $\Sigma$  and stabilize them~\cite{Bro,SS,Suss}. The lack of continuity in the feedback control creates quite a big theoretical problem, because continuity of $U(x)$ is a minimal requirement to apply the classical existence theory of ordinary differential equations to the resulting closed loop system 
\begin{equation}\label{eq:system_cl}
\dot x = f(x,U(x))\,.
\end{equation} 
Therefore, in cases where discontinuous feedback laws have to be used, one has either to choose a generalized concept of solution or to verify that classical solutions still exist when a certain discontinuous law is used.

\smallskip

\n In order to precisely state our results, we need to first introduce a few definitions and notations. Namely, we denote with $|x|$ the Euclidean norm of any element $x\in\R^d$ and with
$$
B_d\doteq\{x\in\R^d~;~|x|< 1\}
$$
the open unit ball of $\R^d$. Also, given any set $E\subseteq \R^d$, we denote its convex hull with $\ch{E}$, and its (topological) closure, interior and boundary respectively with 
$\overline{E}$, $\overset{\,\,\circ}{E}$ and $\partial E$,
so that e.g. we have $\overset{\,\,\circ}{B}_d=B_d$, $\partial B_d=\{x\in\R^d~;~|x|= 1\}$ and $\overline{B_d}=B_d\cup\partial B_d$. 

\smallskip

\n Given a feedback control $u(x)$, we recall that for a system of differential equations like~\eqref{eq:system} with initial datum $x(0)=x_o$, a Carath\'eodory solution on some interval $I$ containing $0$ is an absolutely continuous map $t\mapsto x(t)$ which satisfies \eqref{eq:system} for a.e. $t\in I$, i.e. satisfying the integral representation
\begin{equation}\label{eq:caratheodory}
x(t)=x_0+\int_0^t f(x(s),u(x(s)))\, ds \qquad\forall~t\in I\,.
\end{equation}

\begin{definition}\label{def:ol_control} Given a bounded constraint set $S\subset \R^d$ and a target set $\Sigma\subset \R^d$ such that $S\cap\Sigma\neq\emptyset$, we say that the system~\eqref{eq:system} is \emph{open loop $S$--constrained controllable} to $\Sigma$ if the following holds. For any initial state $x_o\in S$, there exists a Lebesgue measurable control function $u(\cdot)$ and a time $T=T(x_o, u)\geq 0$ such that
denoting with $x(\cdot)$ the Carath\'eodory solution, corresponding to the control $u$, of the Cauchy problem for~\eqref{eq:system} with initial datum $x(0)=x_o$, one has
$$
x(t)\in S\qquad\qquad\forall~t\in [0,T]\,,
$$
and
$$
x(T)\in \Sigma\,.
$$
\end{definition}

\begin{definition}\label{def:cl_stabil} Given a bounded constraint set $S\subset \R^d$ and a target set $\Sigma\subset \R^d$ such that $S\cap\Sigma\neq\emptyset$, we say that a feedback control $U\colon\mathrm{dom}\,U\to{\bf U}$, defined on some open domain $\mathrm{dom}\,U$ which contains $S\setminus \Sigma$, is \emph{$S$--constrained stabilizing} to $\Sigma$ in Carath\'eodory sense for the system~\eqref{eq:system} if the following holds.
For any initial state $x_o\in S$, the closed loop system~\eqref{eq:system_cl} with initial datum $x(0)=x_o$ admits Carath\'eodory solutions and, moreover, for any Carath\'eodory solution $x(\cdot)$ to~\eqref{eq:system_cl} starting from $x_o$ there exists $T\geq 0$ such that one has
$$
x(t)\in S\qquad\qquad\forall~t\in [0,T]\,,\qquad\qquad x(T)\in\Sigma\,.
$$
\end{definition}

One has to be careful when dealing with Definition~\ref{def:cl_stabil}. Indeed, as we have already stressed, in general there might fail to exist a continuous feedback law $U(x)$ which stabilizes~\eqref{eq:system}. Hence, one has to consider discontinuous feedback controls, but in such a case there might be no Carath\'eodory solutions at all.

\n To cope with this problem, we choose here to consider a particular class of feedback controls, the so called \emph{patchy feedbacks}~\cite{AB1, AB2, AB5, BP2}, which are piecewise constant and such that the resulting control system~\eqref{eq:system_cl} always admits local Carath\'eodory solutions for positive times.

\n We remark that this is not the only possible way to overcome the aforementioned difficulties. A possible alternative is to consider time--dependent feedback controls as in~\cite{Cor,Cor2} or to reformulate the problem in the context of hybrid systems~\cite{GPT, GT}. Another approach, closer to the one we adopted here, consists in allowing for arbitrary discontinuous feedback controls $u=u(x)$ and replacing Carath\'eodory trajectories with a weaker concept of solutions. In recent years many authors have followed this alternative path by considering sample--and--hold solutions and Euler solutions for discontinuous vector fields (see e.g.~\cite{CLSW} and references therein) and several results have been obtained in the context of constrained dynamics too (see~\cite{CRS, CS1, CS2}). 

\n However, patchy feedbacks offer some advantage thanks to their robustness. Indeed, one of the main practical problems when using discontinuous controls is the possible appearance of chattering phenomena, which greatly degrade the performance of the control (see e.g.~\cite{Sontag}). A common solution to such problem is the introduction of some observers to filter out the undesired oscillations around the discontinuities~\cite{BlochDrak,DrakUt}. In the context of sample--and--hold solutions for generic measurable feedbacks, a sufficient robustness can be ensured by imposing additional assumptions on the sampling step of the solutions. A nice feature of patchy feedbacks is that none of the above is necessary, because their regularity ensures that only ``tame'' discontinuities are present in the dynamics and thus robustness is guaranteed, provided perturbations are small in the appropriate way (see Theorem~9.4 in~\cite{BPicc} or Theorem~\ref{thm:stab_constr_robust} below). This is in our opinion the main feature of this class of controls: they have a relatively simple structure, being piecewise constant with discontinuities located on the boundaries of a locally finite covering of the state space, they always provide the existence of Carath\'eodory solutions to the dynamics and they offer \emph{complete robustness} w.r.t. both internal and external perturbations, without requiring modifications of the original dynamics and without additional assumptions on their construction procedure.

\smallskip

\n Let us introduce now the assumptions on the constraint set $S$. First we recall a notion from non--smooth analysis~\cite{CLSW}: given a closed set $S$ in $\R^d$ and a point $x\in S$, the \emph{Clarke (proximal) normal cone} to $S$ in $x$ is defined as the set
\begin{equation}\label{eq:normal_cone}
N_S^C(x)\doteq\left\{\lambda\,\xi~;~\lambda\geq 0\,,~\xi\in\chc{\{0\}\cup\Big\{v=\lim_{v_i\to 0}\,{v_i\over |v_i|}\,~;~v_i\perp S \mbox{ in }x_i\,, ~x_i\to x\Big\}}\right\}\,,
\end{equation}
where ``$v_i\perp S$ in $x_i$'' means that $ v_i+x_i\notin S$ and $x_i$ belongs to the projection of $v_i+x_i$ on $S$, or equivalently
$$
v_i+x_i\notin S~\mbox{ and } ~|v_i|=\inf_{\xi\in S} |v_i+x_i-\xi|\,.
$$
We are now ready to state the main hypotheses on $S$:
\begin{description}
\item{{\bf (S1)}} $S$ is compact and wedged at each $x\in\partial S$. The latter means that at each boundary point $x$ one has that $N_S^C(x)$, the Clarke normal cone to $S$ in $x$, is pointed; that is, 
$$
N_S^C(x)\cap\{-N_S^C(x)\}=\{0\}\,.
$$
\item{{\bf (S2)}} The following ``strict inwardness'' condition holds:
$$
\min_{u\in{\bf U}}~ f(x,u)\cdot p~<0\,,
$$
for all $x\in\partial S$ and $p\in N_S^C(x)\setminus\{0\}$.
\end{description}
The former condition is of geometric nature and guarantees that the boundary $\partial S$ is locally homeomorphic to the epigraph of a Lipschitz function~\cite{Rock}. The latter condition states that there always exists at least one control value steering trajectories from the boundary $\partial S$ towards the interior of the constraint set $S$.

\smallskip

Notice that, in Definition~\ref{def:ol_control}, we do not assume \emph{stability} of the target set, i.e. we are not requiring that trajectories starting sufficiently close to $\Sigma$ always remain close to $\Sigma$. Hence, in general we do not expect a feedback which stabilizes the dynamics precisely to the target. The main result of this paper concerns instead \emph{practical stabilization} of~\eqref{eq:system} to $\Sigma$, i.e. the existence for all $\delta>0$ of a patchy feedback control which stabilizes trajectories of~\eqref{eq:system} to a neighborhood $\Sigma^\delta\doteq\Sigma+\delta B_d$ of the target set. 

\begin{theorem}\label{thm:practical_stab_constr} Assume that the system~\eqref{eq:system} satisfies {\bf (F1)}--{\bf (F3)} and open loop $S$--constrained controllability to $\Sigma$, where $S$ is a set satisfying {\bf (S1)} and {\bf (S2)} and $\Sigma$ is any closed set such that $S\cap \Sigma\neq\emptyset$.  Then, for every $\delta>0$ there exists a patchy feedback control $U=U(x)$, defined on an open domain ${\cal D}$ with $S\setminus \Sigma^\delta\subseteq {\cal D}$, which is $S$--constrained stabilizing to $\Sigma^\delta$ for~\eqref{eq:system}.
\end{theorem}

Under assumptions {\bf (S1)} and {\bf (S2)} on the constraint set, it was proved in~\cite{CS1} that it is possible to construct a discontinuous feedback control which steers Euler solutions of~\eqref{eq:system} to $\Sigma^\delta$. 
However, as mentioned above, when a patchy feedback exists more robustness properties can be expected to hold than in the case of the generic feedback presented in~\cite{CS1}. This indeed happens also for the constrained problem considered here, and we will show that practical stabilization of perturbed systems can be established as well. We refer to Section~\ref{sec:rob} for the result about perturbed systems and for further discussions about robustness of the control provided by Theorem~\ref{thm:practical_stab_constr}.

Finally, we want to comment about the convexity assumption {\bf (F3)} on the velocity sets. A large part of our construction does not require it: For instance, only {\bf (F1)} and {\bf (F2)} are needed in Lemma~\ref{lem:patchy_from_wedged} to construct a patchy feedback close to the boundary $\partial S$ of the constraint set.  The only use we make of {\bf (F3)} is to apply the $S$--constrained trajectory tracking lemma from~\cite{CRS} (see Lemma~\ref{lem:track} below), so to replace any $S$--constrained stabilizing open loop control with another control whose corresponding trajectory is stabilized while remaining in the interior of $S$. 
Since we are not aware of an analogous result proved in the case of a general (non--convex) constraint set satisfying {\bf (S1)} and {\bf (S2)} only and of a dynamics satisfying {\bf (F1)} and {\bf (F2)}, but not {\bf (F3)}, this assumption seems unavoidable in Theorem~\ref{thm:practical_stab_constr}. In problems where {\bf (F3)} fails, e.g. when only a discrete range of controls is available, one could reframe the present results in the context of relaxed controls, or could impose additional requirements on the constraint set $S$. In this way, it would be possible to replace Lemma~\ref{lem:track} with a different tracking lemma (without {\bf (F3)}), and the rest of our construction would still hold without further changes.

\section{Preliminiaries}\label{sec:prelim}

In this section we collect a series of preliminary concepts and lemmas which will be used to prove Theorem~\ref{thm:practical_stab_constr}. Namely, in section~\ref{subsec:patchy} we recall the definition of patchy vector fields and patchy feedbacks, and a result (Lemma~\ref{lem:patchy_tubes}) which allows to construct a patchy feedback control around a given trajectory of system~\eqref{eq:system}. Section~\ref{subsec:wedged} contains the definition and the main geometrical properties of sets satisfying {\bf (S1)} and {\bf (S2)}, while section~\ref{subsec:tracking} contains the statement of a ``tracking lemma'' (Lemma~\ref{lem:track}), introduced by Clarke et al. in~\cite{CRS}, which allows to replace $S$--constrained trajectories to~\eqref{eq:system} with other trajectories which always remain in the interior of $S$ and stay close to the original ones. 
Finally, in section~\ref{subsec:patchy_wedged}, we prove a result (Lemma~\ref{lem:patchy_from_wedged}) which allows to construct a patchy feedback close to $\partial S$.

\subsection{Patchy vector fields and patchy feedbacks}\label{subsec:patchy}

We start by recalling the main definitions and properties of the class of discontinuous
vector fields ({\it patchy vector fields}) introduced in~\cite{AB1}.

\begin{figure}
\centering
\psfrag{n}{$\scriptstyle \!\!\!\!\!\!\!{\bf n}(x)$}
\psfrag{g1}{$\scriptstyle \!\!\!\!g_1(x)$}
\psfrag{O1}{$\Omega_1$}
\psfrag{O2}{$\Omega_2$}
\psfrag{O3}{$\Omega_3$}
\includegraphics[width=0.55\textwidth]{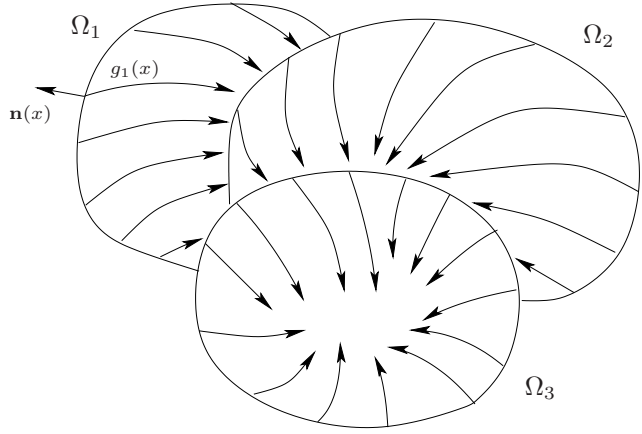}
\caption{
A patchy vector field.}
\label{fig:1}       
\end{figure}

\begin{definition}\label{defn:patchy_vf}
We say that $g:\Omega\to\R^d$ is a \emph{patchy vector field} on the open domain $\Omega\subseteq\R^d$ if there exists a family $\big\{ (\Omega_\alpha,g_\alpha)~;~ \alpha\in\mathcal{A}\big\}$ such that (see Figure~\ref{fig:1})
\begin{description}
\item{(i)} $\mathcal{A}$ is a totally ordered set of indices;
\item{(ii)} each $\Omega_\alpha$ is an open domain with smooth boundary;
\item{(iii)} the open sets $\Omega_\alpha$ form a locally finite covering of $\Omega$;
\item{(iv)} each $g_\alpha$ is a  Lipschitz continuous vector field
defined on a neighborhood of $\overline\Omega_\alpha$, which points strictly inward at each
boundary point $x\in\partial\Omega_\alpha$: namely, calling ${\bf n}(x)$ the outer normal at the boundary
point $x$, we require
\begin{equation}\label{eq:ipc}
g_\alpha(x)\cdot{\bf n}(x) <0\qquad\forall x\in\partial\Omega_\alpha\,;
\end{equation}
\item{(v)} the vector field $g$ can be written in the form
\begin{equation}\label{eq:disc_vf}
g(x) = g_\alpha (x)\qquad \hbox{if}\qquad x \in \Omega_\alpha \setminus \displaystyle{\bigcup_{\beta > \alpha} \Omega_\beta}\,.
\end{equation}
\end{description}
Each element $(\Omega_\alpha,g_\alpha)$ of the family is called \emph{patch}.
\end{definition}

\n By defining
\begin{equation}\label{eq:alpha_star}
\alpha^*(x) \doteq \max\big\{\alpha \in \mathcal{A}~;~ x \in \Omega_\alpha \big\},
\end{equation}
the identity~\eqref{eq:disc_vf} can be written in the equivalent form
\begin{equation}\label{eq:disc_vf_star}
g(x) = g_{_{\alpha^*(x)}}(x) \qquad \forall~x \in \Omega\,.
\end{equation}
We shall occasionally adopt the longer notation
$\big(\Omega,\ g,\ (\Omega_\alpha,\,g_\alpha)_{_{\alpha\in \mathcal{A}}} \big)$
to indicate a patchy vector field, specifying both the domain
and the single patches.

\begin{remark}\label{rem:nonsmooth_bdry} Notice that the smoothness assumption on the boundaries $\partial\Omega_\alpha$ in {\it (ii)} above can be relaxed. Indeed, one can consider patches $(\Omega_\alpha, g_\alpha)$
where the domain $\Omega_\alpha$ only has piecewise smooth boundary. In this case, the inward--pointing condition~\eqref{eq:ipc} can be rephrased as
\begin{equation}\label{eq:ipc_cone}
g(x)\in \overset{\circ}{T}_{\!\Omega}(x)\,,
\end{equation}
$T_{\!\Omega}(x)$ denoting the (Bouligand) tangent  cone to  $\Omega$  at the point $x$, defined by (see~\cite{CLSW})
\begin{equation}\label{eq:Bcone}
T_\Omega(x)\doteq \bigg\{ v\in \R^d~;~ \liminf_{t \downarrow 0} {d\big(x+tv,\ \Omega\big)\over t}=0 \bigg\}\,.
\end{equation}
Clearly, at any regular point $x\in\partial\Omega$, the interior of the tangent cone $T_\Omega(x)$ is precisely the set of all vectors $v\in \R^d$ that satisfy $v\cdot\mathbf{n}(x) <0$ and hence~\eqref{eq:ipc_cone} coincides with the inward--pointing  condition~\eqref{eq:ipc}. 
\end{remark}

\begin{remark} \label{rem:reduced_ipc} Notice also that in Definition~\ref{defn:patchy_vf} the values attained by $g_\alpha$ on $\Omega_\alpha\cap\Omega_\beta$ for any $\beta>\alpha$ are irrelevant. Similarly, the inward pointing condition~\eqref{eq:ipc} does not really matter in points $x\in\partial\Omega_\alpha\cap\Omega_\beta$, for any $\beta>\alpha$, and in points $x\in\partial\Omega_\alpha\cap (\R^d\setminus \Omega)$. This is a consequence of the fact that, in general, the patches $(\Omega_\alpha,\,g_\alpha)$ are not uniquely determined by the patchy vector field $g$.

\n Indeed, as observed in~\cite{AB1}, whenever a Lipschitz vector field $h_\alpha$ is given on $\overline \Omega_\alpha$ so that it verifies~\eqref{eq:ipc} on $(\partial \Omega_\alpha\cap\Omega)\setminus \bigcup_{\beta>\alpha}\Omega_\beta$, one can always construct another Lipschitz vector field $g_\alpha$ on $\overline \Omega_\alpha$ such that $g_\alpha=h_\alpha$ on $(\overline\Omega_\alpha\cap\Omega)\setminus \bigcup_{\beta>\alpha}\Omega_\beta$ and such that~\eqref{eq:ipc} is verified at every $x\in\partial\Omega_\alpha$.
\end{remark}

If $g$ is a patchy vector field, the differential equation
\begin{equation}\label{eq:patchy_system}
\dot x = g(x)
\end{equation}
has many interesting properties.
In particular, it was proved in~\cite{AB1} that, given any initial condition
\begin{equation}\label{eq:patchy_datum}
x(0)=x_0,
\end{equation}
the Cauchy problem \eqref{eq:patchy_system}--\eqref{eq:patchy_datum} has at least one forward
solution, and at most one backward solution in Carath\'eodory sense. We recall that a Carath\'eodory solution of~\eqref{eq:patchy_system}--\eqref{eq:patchy_datum} on some interval $I$ is an absolutely continuous map $t\mapsto\gamma(t)$ which satisfies \eqref{eq:patchy_system} for a.e. $t\in I$, i.e.
\begin{equation}\label{eq:caratheodory_patchy}
\gamma(t)=x_0+\int_0^t g(\gamma(s))\, ds \qquad\forall~t\in I\,.
\end{equation}

\n We collect below the other main properties satisfied by trajectories of~\eqref{eq:patchy_system}--\eqref{eq:patchy_datum}.
\begin{itemize}
\item For every Carath\'eodory solution $\gamma(\cdot)$ of~\eqref{eq:patchy_system}, the map $t \mapsto \alpha^*(\gamma(t))$, with $\alpha^*$ the function defined in~\eqref{eq:alpha_star}, is left continuous and non-decreasing. Moreover, it is piecewise constant on every compact interval $[a,b]$, i.e. there exist a partition $a=t_o<t_1<\ldots<t_N=b$ of $[a,b]$ and indices $\alpha_1<\ldots<\alpha_N$ in ${\cal A}$ such that $\alpha^*(\gamma(t))=\alpha_i$ for all $t\in\,]t_{i-1},t_i]$.
\item The set of all Carath\'eodory solutions of~\eqref{eq:patchy_system}--\eqref{eq:patchy_datum} is closed in the topology of uniform convergence, but possibly not connected.
\item Carath\'eodory solution of~\eqref{eq:patchy_system} are robust w.r.t. to both inner and outer perturbations; namely, for any solution $y(\cdot)$ of the perturbed system
$$
\dot y=g(y+\zeta)+d\,,
$$
there exists a solution $x(\cdot)$ of the unperturbed system~\eqref{eq:system} such that $||x-y||_{{\bf L}^\infty}$ is as small as we want, provided that $\zeta$ and $d$ are small enough in $\BV$ and $\Lsp^1$, respectively (see~\cite{AB2} for the details in the general case, and Section~\ref{sec:final} for a discussion of the constrained case).
\end{itemize}

\medskip

The class of patchy vector fields is of great interest in a wide variety of control problems for  general nonlinear control systems~\eqref{eq:system}, that can be solved by constructing a state feedback $u=U(x)$ which renders the resulting closed loop map $g(x)=f(x, U(x))$ a patchy vector field and, hence, ensures robustness properties of the resulting solutions without additional efforts. This leads to the following definition.
\begin{definition}\label{defn:patchy_feedback}
Let $\mathcal{A}$ be a totally ordered set of indices, and $(U_\alpha)_{\alpha\in\mathcal{A}}$ be a family of control values in ${\bf U}$ such that, for each $\alpha\in\mathcal{A}$, there exists a patch $(\Omega_\alpha, g_\alpha)$ which satisfies
\begin{equation}\label{eq:patchy_feedback_1}
g_\alpha(x) = f(x,\, U_\alpha)\qquad\qquad\forall~x \in
\Omega_\alpha \setminus \bigcup_{\beta > \alpha} \Omega_\beta\,.
\end{equation}
If the family $\{\Omega_\alpha\}_{\alpha\in\mathcal{A}}$ forms a locally finite covering of an open domain ${\cal D}\subseteq\R^d$, then the piecewise constant map
\begin{equation}\label{eq:patchy_feedback_2}
U(x) \doteq  U_\alpha\qquad \hbox{if}\qquad x \in
\Omega_\alpha \setminus \bigcup_{\beta > \alpha} \Omega_\beta
\end{equation}
is called a \emph{patchy feedback control} on ${\cal D}$.
\end{definition}

By requiring~\eqref{eq:patchy_feedback_1} with $(\Omega_\alpha, g_\alpha)$ being a patch, in particular we require that
$$
f(x,\, U_\alpha(x))\cdot{\bf n}(x) <0\qquad\forall x\in\partial\Omega_\alpha \setminus \bigcup_{\beta > \alpha} \Omega_\beta\,.
$$
By Definitions~\ref{defn:patchy_vf}--\ref{defn:patchy_feedback} it is thus clear that, given
a patchy feedback $U$, the corresponding collection of patches $(\Omega_\alpha, g_\alpha)$, $\alpha\in\mathcal{A}$, defines a patchy vector field $g(x)=f(x, U(x))$ on $\bigcup_{\alpha\in\mathcal{A}}\Omega_\alpha$.
Moreover, recalling the definition of $\alpha^*(x)$ in \eqref{eq:alpha_star}, a patchy feedback control can be written in the
equivalent form
\begin{equation}\label{eq:patchy_feedback_2_star}
U(x)=U_{\alpha^*(x)}(x)\qquad\qquad x\in \Omega\doteq\bigcup_{\alpha\in\mathcal{A}}\Omega_\alpha\,.
\end{equation}
We shall occasionally adopt the longer notation $(U,\, (\Omega_\alpha, U_\alpha)_{\alpha\in\mathcal{A}})$ to indicate a patchy feedback control, similarly to the notation adopted for patchy vector fields.

\begin{remark}\label{rem:patchy_nonconst}
As in Remark~\ref{rem:reduced_ipc}, the values attained by $U_\alpha$ on the set $\Omega_\alpha\cap\Omega_\beta$ are irrelevant, whenever $\alpha<\beta$, and similarly it only matters that $f(\cdot,U_\alpha(\cdot))$ fulfills the inward--pointing condition~\eqref{eq:ipc}  at points of $(\partial\Omega_\alpha\cap\Omega)\setminus \bigcup_{\beta>\alpha}\Omega_\beta$.
\end{remark}

We now recall a result, originally proved in~\cite{AB1}, which allows to construct a patchy feedback control, starting from a piecewise constant open loop control.


\begin{lemma}\label{lem:patchy_tubes} Assume that $f$ in~\eqref{eq:system} satisfies {\bf (F1)} and {\bf (F2)}. Let $x_o\in \R^d$, $T>0$ and $u\colon [0,T]\to{\bf U}$ be a piecewise constant open loop control. Then, denoting by $x(\cdot)$ the solution to~\eqref{eq:system} with initial datum $x_o$ and control $u(\cdot)$, for every $\ve>0$ there exists an open domain $\Gamma\subset \R^d$ and a patchy feedback control $U$ defined on a domain ${\cal D}\supseteq\Gamma$ such that the following properties hold.
\begin{description}
\item{\it (i)} There holds
$$
\bigcup_{t\in [0,T]}x(t)~\subseteq~ \Gamma~\subseteq~\bigcup_{t\in [0,T]}x(t)+\ve B_d\,.
$$
\item{\it (ii)} The patchy feedback $U$ coincides with $u$ along the trajectory $x(\cdot)$, namely $U(x(t))=u(t)$ for a.e. $t\in [0,T]$.
\item{\it (iii)} The patchy vector field $\xi\mapsto g(\xi)\doteq f(\xi,U(\xi))$ satisfies
$$
g(z)\in \overset{\circ}{T}_{\!\Gamma}(z)
\qquad\forall z\in\partial\Gamma\setminus B^+\,,
$$
where $T_{\!\Gamma}(z)$ denotes the Bouligand tangent cone to $\Gamma$ at the point $z$, as in~\eqref{eq:Bcone},
and 
$$
B^+\doteq\left\{\eta\in x(T)+\ve B_d~;~(\eta-x(T))\cdot \dot x(T)\geq 0\right\}\,.
$$
In other words, $g$ satisfies the inward--pointing condition in all points of $\partial\Gamma$, except possibly those which lie in a half--ball of radius $\ve$ centered at the end point $x(T)$ of the trajectory.
\item{\it (iv)} For all $\xi\in\overline{\Gamma}$ and all Carath\'eodory solutions $y(\cdot)$ to
$$
\dot y(t)=g(y(t))\,\qquad\qquad y(0)=\xi\,,
$$
there exists $\tau\in[0,T]$ such that 
$y(\cdot)$ is defined on $[0,T-\tau]$, $y(t)\in\Gamma$ for all $t\in \,]0,T-\tau]$ and there holds
$$
\max_{t\in [0,T-\tau]}|x(\tau+t)-y(t)|\leq \ve\,.
$$
In particular, all trajectories starting in $\Gamma$ and corresponding to the patchy feedback $U$ remain close to $x(\cdot)$ and eventually reach $x(T)+\ve\overline{B_d}$.
\end{description}
\end{lemma}

\begin{figure}
\centering
\psfrag{x}{$\scriptscriptstyle \!\!\!\!x(0)$}
\psfrag{g}{$\scriptscriptstyle \!\!\!\!x(t)$}
\psfrag{u}{$\scriptscriptstyle x(T)$}
\psfrag{t}{$\scriptstyle \Gamma_o$}
\psfrag{v}{$\scriptstyle \Gamma_1$}
\includegraphics[width=0.65\textwidth]{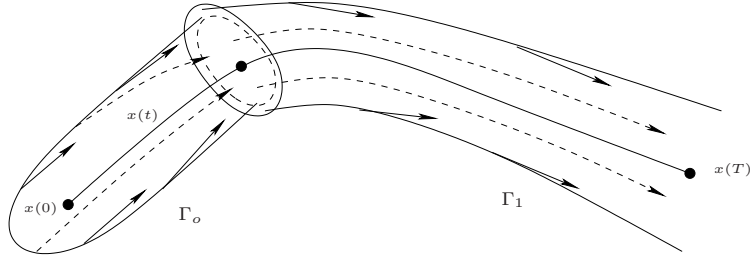}
\caption{Tube--like domain $\Gamma$ around a solution $x(\cdot)$ with piecewise constant open loop control.}
\label{fig:tubes}       
\end{figure}

The basic idea is that, once you are given a piecewise constant control $u(\cdot)$ and the corresponding trajectory
$x(\cdot)$, if you use the same control values in a small neighborhood of the trajectory then the resulting solutions of~\eqref{eq:system} will stay close to  $x(\cdot)$ itself, by regularity of the vector field, see Figure~\ref{fig:tubes}. Hence, one can collect in a set $\Gamma_o$ all points of trajectories obtained with constant control $u(0)$ and initial datum close to $x_o$, say in $x_o+\rho_oB_d$, until the first time $t_1$ in which $u(\cdot)$ has a jump; then restart the procedure and collect in a set $\Gamma_1$ all points of trajectories obtained with constant control $u(t_1)$ and initial datum in $x(t_1)+\rho_1B_d$ until next time $t_2$ in which $u(\cdot)$ has a jump; and so on. Since the number of jumps in $u$ is finite, if the radii $\rho_j>0$ are chosen properly, then the finite union of reachable sets $\Gamma=\bigcup_{j=0}^N\Gamma_j$ has the required properties. The detailed proof of this lemma can be found as part of the proof of Proposition~4.1 in~\cite{AB1}, or as part of the proof of Lemma~9.1 in~\cite{BPicc}.

\subsection{Wedged sets and inner approximations}\label{subsec:wedged}

In this section we collect a few geometrical properties which are satisfied by either any wedged set or specifically by sets for which both {\bf (S1)} and {\bf (S2)} hold. 

As a first step, we want to give a characterization of wedged sets. We start by recalling the definition of Clarke's tangent cone to $S$ in $x$, which we denote by $T_S^C(x)$, as the polar cone to the normal cone $N_S^C(x)$ defined in~\eqref{eq:normal_cone}, that is
\begin{equation}\label{eq:Ccone}
T_S^C(x)\doteq\left\{v\in\R^d~;~p\cdot v\leq 0~~~\forall p\in N_S^C(x)\right\}\,.
\end{equation}
In general, the Clarke tangent cone is smaller than the Bouligand tangent cone defined in~\eqref{eq:Bcone} (see~\cite{CLSW}), 
i.e. there holds for every closed set $S$ and every $z\in S$
$$
T^C_S(z)~\subseteq~ T_S(z)\,.
$$

\begin{figure}
\centering
\psfrag{e}{$\scriptstyle \!\!\!\!\ve |v|$}
\psfrag{e2}{$\scriptstyle \ve$}
\psfrag{v}{$\scriptstyle v$}
\psfrag{B}{$\scriptstyle v+\ve B_d$}
\psfrag{W}{$\scriptstyle {\cal W}(v,\ve)$}
\psfrag{w1}{$\scriptstyle {\cal W}(v,\ve)$}
\psfrag{w2}{$\scriptstyle {\cal W}(v,\ve/2)$}
\psfrag{P}{$\scriptstyle \!\!\!\!\!\!\!\!\partial^-{\cal W}(v,\ve/2)$}
\includegraphics[width=0.65\textwidth]{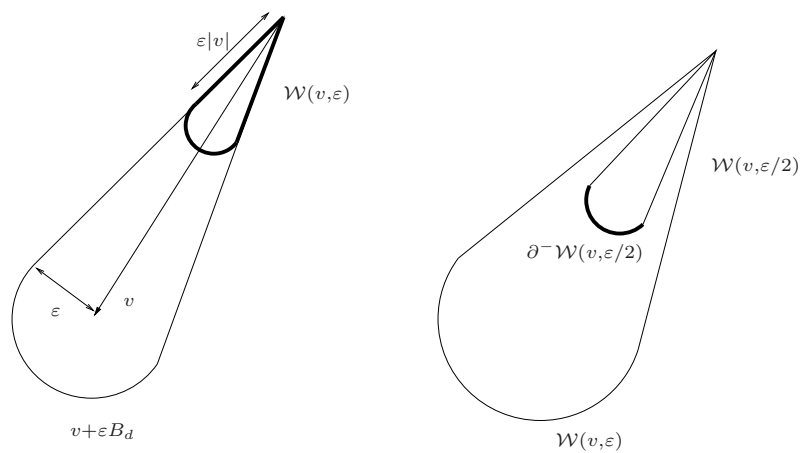}
\caption{
Left: Wedge ${\cal W}(v,\ve)$, of axis $v$ and radius $\ve$. Right: The ``lower'' boundary $\partial^-{\cal W}(v,\ve/2)$ of the smaller wedge is strictly separated from $\R^d\setminus{\cal W}(v,\ve)$.}
\label{fig:6}
\end{figure}

\n Moreover, for all $v\in\R^d$ and $\ve>0$ we call \emph{wedge} of axis $v$ and radius $\ve$ the set (see Figure~\ref{fig:6} left)
$$
{\cal W}(v,\ve)\doteq\left\{sw~;~w\in v+\ve B_d,\,s\in [0,\ve]\right\}\,.
$$
Finally, to denote  the ``lower'' part of the boundary of a wedge (see Figure~\ref{fig:6} right), we use the following
$$
\partial^-{\cal W}(v,\ve)\doteq\left\{\ve w~;~w\in v+\ve \partial B_d,\, (v-w)\cdot v\leq 0\right\}\,.
$$

\smallskip

We are now in a position to state the following characterization result.

\begin{proposition}\label{prop:wedged} Let $S\subseteq\R^d$ be a closed nonempty set and $x\in\partial S$. Then, the following properties are equivalent:
\begin{description}
\item{{\it (i)}} $S$ is wedged in $x$;
\item{{\it (ii)}} $T_S^C(x)$ has nonempty interior;
\item{{\it (iii)}} there exist $v\in\R^d$ and $\ve>0$ such that
$$
y+{\cal W}(v,\ve)\subset S\qquad\qquad \forall~y\in \{x+\ve B_d\}\cap S\,.
$$
\end{description}
\end{proposition}
We refer to~\cite{CLSW} for the proof of the equivalences above. 
We also recall that for a wedged set $S$ the following properties hold for every $x\in\partial S=\partial(\R^d\setminus S)$
\begin{equation}\label{eq:cones_to_wedged}
N_S^C(x)=-N_{\R^d\setminus S}^C(x)\,,
\qquad\qquad \qquad
T_S^C(x)=-T_{\R^d\setminus S}^C(x)\,.
\end{equation}
Finally, we mention that wedged sets are sometimes called \emph{epi--Lipschitz sets} because they are locally the epigraph of a Lipschitz continuous function (see~\cite{Rock}).

\n For later use, we need a better understanding of the behavior of wedges when their radii are rescaled. It is immediate to deduce from the definition that ${\cal W}(v,\ve/2)\subseteq{\cal W}(v,\ve)$ for every $v\in\R^d$ and $\ve>0$. Moreover, we claim that points of the ``lower'' boundary $z\in\partial^-{\cal W}(v,\ve/2)$, are well inside the larger wedge ${\cal W}(v,\ve)$ (see again Figure~\ref{fig:6} right). Indeed, it is not difficult to verify that for any fixed $v\in\R^d$ and any $z\in\partial^-{\cal W}(v,\ve/2)$, 
there holds
\begin{equation}\label{eq:rescaled_wedge}
{\ve\over {\cal O}(1)}\,\leq\, d(z,\overline{\R^d\setminus{\cal W}(v,\ve)})\,\leq\,\ve\,  {\cal O}(1)\,.
\end{equation}

\medskip

Next, we recall a result concerning the ``inner approximations'' of a set $S$ satisfying {\bf (S1)}. Given a closed set $S$ and $r\geq 0$, we call \emph{$r$--inner approximation of $S$} the set
\begin{equation}\label{eq:inner_ap}
S_r\doteq\{x\in\R^d~;~d(x,\R^d\setminus S)\geq r\}\,,
\end{equation}
and we define
\begin{equation}\label{eq:crown}
Q(S,r)\doteq S\setminus \overset{\!\!\circ}{S_r}\,.
\end{equation}
Given $x\in S$, we also set $r(x)\doteq d(x,\overline{\R^d\setminus S})$,
and we introduce the following notations
$$
N(x)\doteq \left\{v \in N^C_{S_{r(x)}}(x)~;~ |v|=1\right\}\,,
\qquad\qquad
T(x)\doteq T^C_{S_{r(x)}}(x)\,.
$$

\begin{lemma}[Lemma~3.3 in~\cite{CRS}]\label{lem:CRS1} Let $S$ be a set such that {\bf (S1)} is verified. Then there exists $r_o>0$ such that for all $r\in[0,r_o]$ the set $S_r$ is nonempty and wedged in every point of its boundary. Moreover, the multifunction $N(\cdot)$ has closed graph on $Q(S,r_o)$ or, equivalently, the multifunction $T(\cdot)$ is lower semicontinuous on $Q(S,r_o)$.
\end{lemma}
 
Assumption {\bf (S2)} on the constraint set plays an important role as well. Observe that in terms of Clarke's tangent cone~\eqref{eq:Ccone} to a closed set $S$, condition {\bf (S2)} can be restated as follows. For every $x\in\partial S$ there holds
\begin{equation}\label{eq:pos_invariance}
\overset{\!\!\!\!\circ}{T^C_S}(x)\cap f(x,{\bf U})\neq \emptyset\,,
\end{equation}
i.e. {\bf (S2)} ensures that there is an admissible speed pointing strictly inside the set $S$. 
\medskip

Finally, we mention a result dealing with decrease properties of the signed distance function. We recall that, given a closed set $Z\subset \R^d$, the \emph{signed distance} of a point $x\in\R^d$ from $Z$ is defined by
\begin{equation}\label{eq:sign_dist}
\Delta_Z(x)\doteq d(x,Z)-d(x,\overline{\R^d\setminus Z})\,.
\end{equation}
It is not difficult to verify that the function $\Delta_Z$ is Lipschitz continuous.

\begin{lemma}[Lemma~3.7 in~\cite{CRS}]\label{lem:CRS4} Let $Z$ be a closed set which is wedged at $x\in\partial Z$. 
Assume $v\in \R^d$ and $\ve>0$ are such that
$$
y+{\cal W}(v, \ve)\subset Z\qquad\qquad \forall~y\in \{x+\ve B_d\}\cap Z\,,
$$
$$
y+{\cal W}(-v, \ve)\subset \overline{\R^d\setminus Z}\qquad\qquad \forall~y\in \{x+\ve B_d\}\setminus \overset{\,\circ}{Z}\,.~~~~~~~~~
$$
Then there exists a neighborhood ${\cal N}_x$ of $x$ such that
$$
\nabla\Delta_Z(y)\cdot v \leq -\ve\,,
$$
for all $y\in{\cal N}_x\setminus\partial S$ in which $\Delta_Z$ is differentiable.
\end{lemma}

\subsection{Trajectory tracking}\label{subsec:tracking}

Since in the following results we need to compare trajectories starting in the same point with different controls, or from different points using the same control, throughout this section we will denote the solution to the Cauchy problem
$$
\dot x = f(x,u)\,,\qquad\qquad\qquad x(t_o)=x_o\,,
$$
with $x(\cdot\,;t_o,x_o,u(\cdot))$, for every $t_o\in\R$, $x_o\in\R^d$ and $u\colon[t_o,t_1]\to{\bf U}$ measurable open loop control. First of all, we recall the \emph{trajectory tracking lemma} introduced in~\cite{CRS}. 

\begin{lemma}[Lemma~3.9 in~\cite{CRS}]\label{lem:track} Assume that $f$ in~\eqref{eq:system} satisfies {\bf (F1)}--{\bf (F3)}. Let $S$ be a set such that {\bf (S1)} and {\bf (S2)} are verified and $r_o>0$ be the value found in Lemma~\ref{lem:CRS1}. Then there exist constants $C>0$ and $T^*>0$ such that for every $r\in [0,r_o]$ the following holds. Given any initial datum $x_o$ in $S_r$ and any measurable open loop control $u(\cdot)$, there exists another measurable open loop control $\bar u(\cdot)$ such that
$$
x(t\,;0,x_o,\bar u(\cdot))\in S_r\qquad\qquad \forall~t\in [0,T^*]\,,
$$
and
$$
\big|x(t\,;0,x_o,u(\cdot))-x(t\,;0,x_o,\bar u(\cdot))\big|~\leq ~C\max_{s\in[0,T^*]} d(x(s\,;0,x_o,u(\cdot)), S_r)\qquad\qquad \forall~t\in [0,T^*]\,.
$$
\end{lemma}

In the present work, we use the tracking lemma only to prove the following approximation result for $S$--constrained trajectories by means of $S_r$--constrained trajectory, for a sufficiently small $r>0$.
Its proof is contained in the one of Theorem~3.10 in~\cite{CRS}.

\begin{corollary}\label{cor:stick_inside} Assume that $f$ in~\eqref{eq:system} satisfies {\bf (F1)}--{\bf (F3)}. Let $S$ be a set such that {\bf (S1)} and {\bf (S2)} are verified and $r_o>0$ be the value found in Lemma~\ref{lem:CRS1}. Given $T>0$, there exists a constant $C_T>0$ such that for any initial state $\xi\in\overset{\,\,\circ}{S}$ the following holds. If $r\in[0,r_o]$ is such that $\xi\in S_r$ and $u\colon[0,T]\to{\bf U}$ is a measurable open loop control such that the corresponding trajectory of~\eqref{eq:system} emanating from $\xi$ satisfies
\begin{equation}\label{eq:app_traj_hyp}
x(t\,;0,\xi,u(\cdot))\in S\qquad\qquad\forall~t\in[0,T]\,,
\end{equation}
then there exists a control $\bar u(\cdot)$ defined on $[0,T]$ such that the corresponding solution $x(\cdot\,;0,\xi,\bar u(\cdot))$ satisfies $\forall~t\in[0,T]$
\begin{equation}\label{eq:app_traj_01}
x(t\,;0,\xi,\bar u(\cdot))\in S_r\,,
\end{equation}
\begin{equation}\label{eq:app_traj_02}
\big|x(t\,;0,\xi,\bar u(\cdot))-x(t\,;0,\xi,u(\cdot))\big|~\leq~C_T\,r\,.
\end{equation}
\end{corollary}

\smallskip

\subsection{Patchy feedback close to $\partial S$}\label{subsec:patchy_wedged}

In this section we present a technical lemma dealing with the construction of a patchy feedback control $U(x)$ near the boundary of a wedged set $S$ satisfying {\bf (S1)} and {\bf (S2)}, so that $S$ results positively invariant for the resulting dynamics~\eqref{eq:system_cl}. Notice that the sets $S_r$ and $Q(S,r)$, for any $r>0$, have been introduced in~\eqref{eq:inner_ap}--\eqref{eq:crown} and that for any $x\in S$ we still use $r(x)$ as a shorter notation for the quantity $d(x,\overline{\R^d\setminus S})$. Also notice that the following result does not use assumption {\bf (F3)}.

\begin{lemma}\label{lem:patchy_from_wedged} Assume that $f$ in~\eqref{eq:system} satisfies {\bf (F1)} and {\bf (F2)} and that $S$ satisfies {\bf (S1)} and {\bf (S2)}. Let $r_o>0$ be as in Lemma~\ref{lem:CRS1}.
Then, there exist $\tilde r\in\,]0,r_o[$, $\mu>0$ and a patchy feedback control $U\colon {\cal D} \to {\bf U}$, with 
\begin{equation}\label{eq:domain_patchy_wedged}
Q(S, \tilde r)\subseteq  {\cal D}\,,
\end{equation}
such that for all $x\in Q(S, \tilde r)$ there holds
\begin{equation}\label{eq:inward_bdry}
f(x,U(x))+\mu\, B_d~\subset~ T^C_{S_{r(x)}}(x)\,.
\end{equation}
Moreover, one can require that for every $\lambda>0$ the patchy control $U=(U,\, (\Omega_k, U_k)_{k\in\mathcal{K}})$ satisfies $\diam\,\Omega_k\leq \lambda$.
\end{lemma}

\medskip

\n{\it Proof.} {\bf Step~1}. Let $x\in\partial S$. Since~\eqref{eq:pos_invariance} holds, there exist $u_x\in{\bf U}$ and $\mu_x>0$ such that
\begin{equation}\label{eq:ipc_constraint}
f(x,u_x)+\mu_x B_d~\subset~ T^C_S(x)=T(x)\,.
\end{equation}
Hence, by setting for ease of notation $w(x)\doteq f(x, u_x)$, there exists $\tilde\ve=\tilde\ve(x)>0$ such that
\begin{equation}\label{eq:int_cone_vertex}
y+{\cal W}\big(w(x), \tilde\ve\big)~\subset ~S~~~~~~~~\qquad\qquad \forall~y\in \{x+\tilde\ve B_d\}\cap S\,,
\end{equation}
Notice that
$$
|w(x)|\leq M\doteq\max\big\{|f(x,u)|~;~(x,u)\in{S\times{\bf U}}\big\}\,,
$$
and that we must have $|w(x)|>0$, because otherwise from~\eqref{eq:ipc_constraint} we would have $T^C_S(x)=\R^d$, and thus $N^C_S(x)=\{0\}$. 

We also deduce that, by continuity of $f$ and by Lemma~\ref{lem:CRS1} (lower semicontinuity of the multifunction $T(\cdot)$ on $Q(S,r_o)$), there exists $\rho_x>0$ such that for all $\xi\in S$ with $|x-\xi|<\rho_x$ there holds
\begin{equation}\label{eq:ipc_constraint2}
f(\xi,u_x)+\,{\mu_x\over 2}\, B_d~\subset~ T(\xi)\,.
\end{equation}
Notice that, by possibly reducing $\rho_x$, it is not restrictive to assume $\rho_x<\lambda/2$, for a fixed constant $\lambda>0$, and to assume that for all $\xi$ with $|x-\xi|<\rho_x$ there also holds $|f(\xi,u_x)-w(x)|<\tilde\ve/4$.

\smallskip

Now fix $\beta=\beta(x)$ such that
$0<\beta<\min\left\{
1\,,\,
{2\,\rho_x\over 2 M \tilde \ve+\tilde\ve^2}
\right\}$,
and consider the rescaled wedge $\beta{\cal W}(w(x),\tilde\ve/2)$.
\begin{figure}
\centering
\psfrag{Q}{$\scriptstyle S$}
\psfrag{S}{$\scriptstyle S_{r_o}$}
\psfrag{x}{$\scriptscriptstyle x$}
\psfrag{w1}{$\scriptscriptstyle \!\!{\cal W}(w(x),\tilde\ve)$}
\psfrag{w}{$\scriptscriptstyle \!\!\beta{\cal W}(w(x),\tilde\ve/2)$}
\psfrag{r1}{$\scriptscriptstyle \!\!\!\!\partial S$}
\psfrag{r2}{$\scriptscriptstyle \!\!\!\!S_{R}$}
\includegraphics[width=0.42\textwidth]{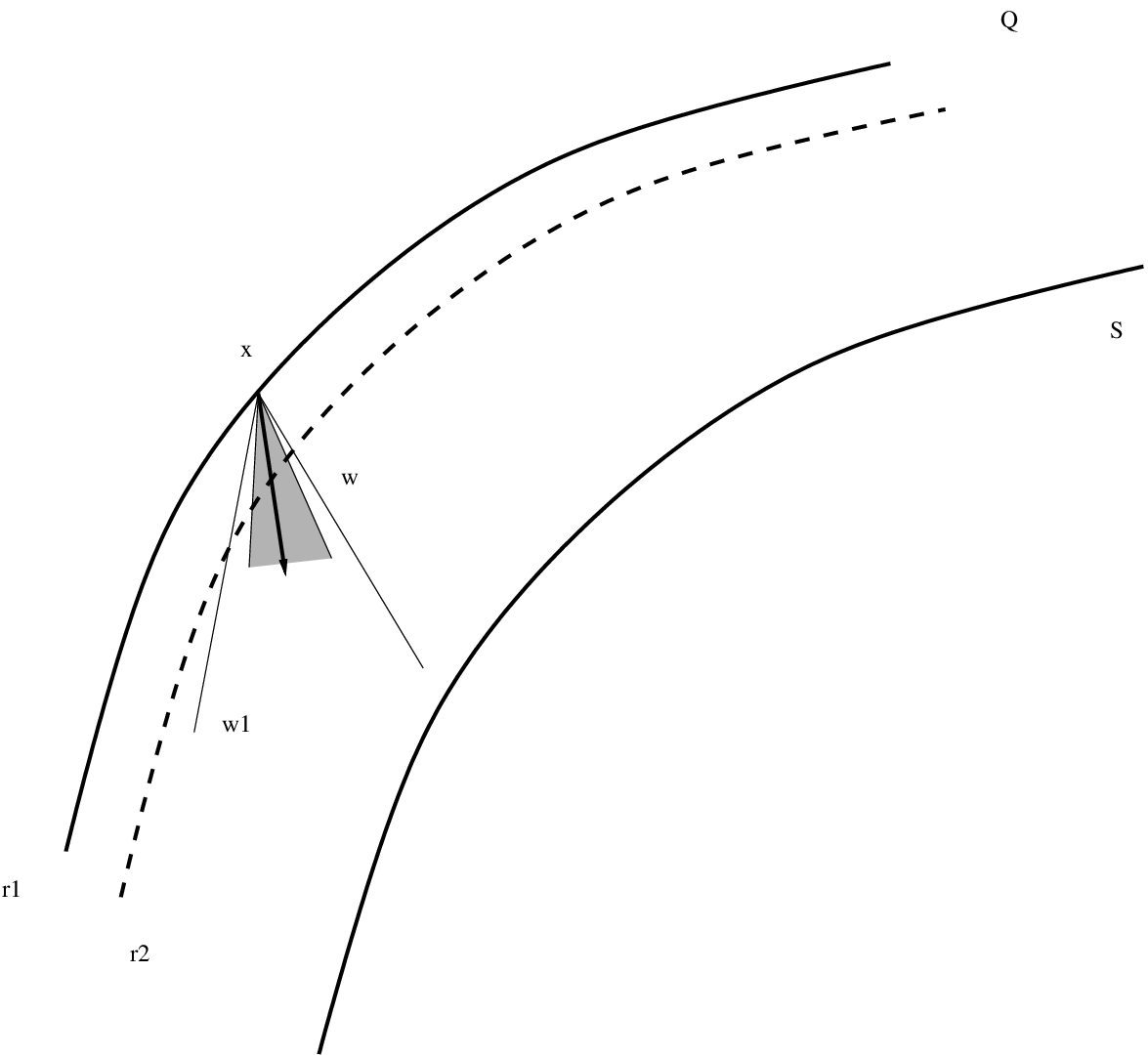}
~~~~~~~~~
\psfrag{Q}{$\scriptstyle S$}
\psfrag{S}{$\scriptstyle S_{r_o}$}
\psfrag{x}{$\scriptscriptstyle x$}
\psfrag{w}{$\scriptscriptstyle \Gamma^x$}
\psfrag{v}{$\scriptscriptstyle w(x)$}
\psfrag{r1}{$\scriptscriptstyle \!\!\!\!\partial S$}
\psfrag{r2}{$\scriptscriptstyle \!\!\!\!S_{\alpha}$}
\includegraphics[width=0.42\textwidth]{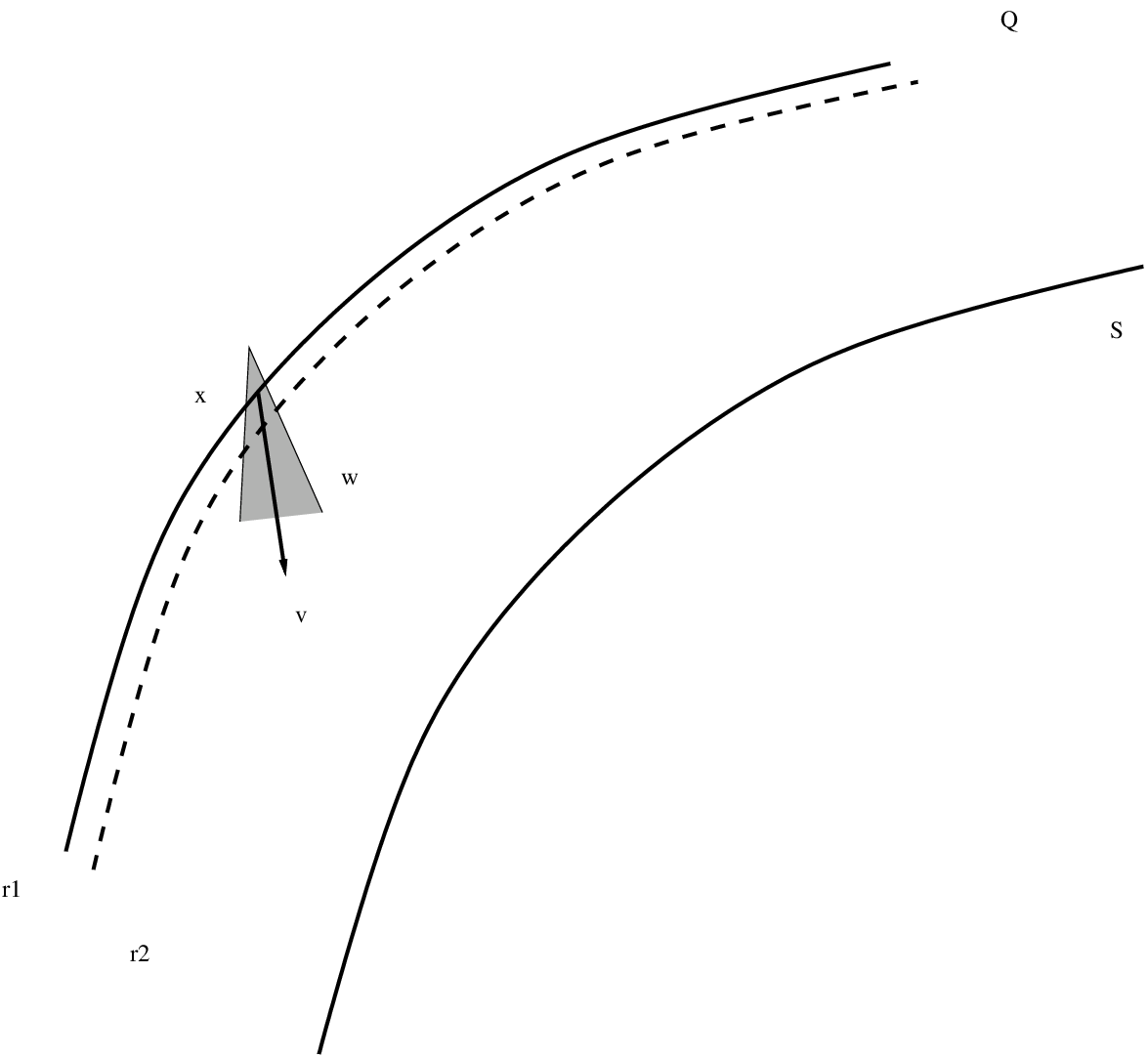}
\caption{
Left: The rescaled wedge $x+\beta{\cal W}(w(x),\tilde\ve/2)$ is well inside the wedge $x+{\cal W}(w(x),\tilde\ve)$. Right: The domain $\Gamma^x$ is obtained by a small translation of $x+\beta{\cal W}(w(x),\tilde\ve/2)$.}
\label{fig:8}
\end{figure}
Recalling that for points $z\in\partial^-{\cal W}(w(x),\tilde\ve/2)$ there holds~\eqref{eq:rescaled_wedge}, we obtain a completely similar estimate for points
$z\in\partial^-\Big(\beta{\cal W}(w(x),\tilde\ve/2)\Big)=\beta \partial^-{\cal W}(w(x),\tilde\ve/2)$.
Namely, every point in $z\in\partial^-\Big(\beta{\cal W}(w(x),\tilde\ve/2)\Big)$ has an uniformly positive distance from $\R^d\setminus{\cal W}(w(x),\tilde\ve)$, that we denote with $R=R(x)$, so that (see Figure~\ref{fig:8} left) 
$$
0<R\leq d(z,\overline{\R^d\setminus{\cal W}(w(x),\tilde\ve)})\qquad\qquad\forall z\in\partial^-\Big(\beta{\cal W}(w(x),\tilde\ve/2)\Big)\,.
$$

We finally fix $\alpha=\alpha(x)$ such that
$0<\alpha<\min\left\{
{\rho_x\over 2}\,,\,
{R\over 2}\,,\,
{\beta\over 2}
\right\}$, and we define
$$
\Gamma^x\doteq x-\alpha\,{w(x)\over |w(x)|}+ \beta \,\overset{\,\,\circ}{\cal W}\left(w(x),\tilde\ve/2\right)\,.
$$
In other words, the domain $\Gamma^x$ is obtained by first rescaling the (open) wedge $x+\overset{\,\,\circ}{\cal W}(w(x),\tilde\ve/2)$ by a suitable factor $\beta$ and then by slightly shifting it along the direction $-w(x)$, so to obtain an open neighborhood of $x$. The particular choice of $\alpha, \beta$ ensures that there holds
\begin{equation}\label{eq:inclusion_lemma2}
\Gamma^x\subseteq x+\rho_x B_d\,,
\end{equation}
so that~\eqref{eq:ipc_constraint2} holds in every point of $\Gamma^x$.

\n If we split the boundary $\partial\Gamma^x$ into its ``lower'' and ``upper'' parts, by setting
$$
\partial^-\Gamma^x\doteq x-\alpha\,{w(x)\over |w(x)|}\, + \partial^-\Big(\beta{\cal W}(w(x),\tilde\ve/2)\Big)\,,
\qquad\qquad
\partial^+\Gamma^x\doteq\partial\Gamma^x\setminus \partial^-\Gamma^x\,,
$$
the choice of the parameters $\alpha$ and $\beta$ also ensures that the lower boundary $\partial^-\Gamma^x$ is well inside the wedge $x+{\cal W}(w(x),\tilde\ve)$. Indeed, for $z\in\partial^-\Gamma^x$ we have
\begin{align*}
d\left(z~,~\overline{\R^d\setminus\big(x+{\cal W}(w(x),\tilde\ve)\big)}\,\right)&\geq d\left(z+\alpha\,{w(x)\over |w(x)|}~,~\overline{\R^d\setminus\big(x+{\cal W}(w(x),\tilde\ve)\big)}\,\right)
-\alpha\\
&>R-\alpha> \alpha>0\,.
\end{align*}
In turn, since the wedge is contained in $S$ by {\bf (S1)}, this implies for every $z\in\partial^-\Gamma^x$
\begin{align*}
d\left(z~,~\overline{\R^d\setminus S}\,\right)&\geq d\left(z~,~\overline{\R^d\setminus\big(x+{\cal W}(w(x),\tilde\ve)\big)}\,\right)>\alpha\,,
\end{align*}
which ensures $z\in S_\alpha$. Hence, by replacing $\Gamma^x$ with
$\Omega^x\doteq\Gamma^x\setminus S_\alpha$,
we still obtain an open neighborhood of $x$, contained into $x+\rho_x B_d$, whose boundary can be divided into
$$
\partial^+\Omega^x\doteq\partial^+\Gamma^x\setminus S_\alpha\,,
\qquad
\qquad
\partial^-\Omega^x\doteq\partial\Omega^x\setminus\partial^+\Omega^x=\Gamma^x\cap\partial S_\alpha\,,
$$
and $r(z)=\alpha$ for all $z\in\partial^-\Omega^x$. 

\medskip

\n{\bf Step~2}. 
Now, by compactness of $\partial S$, we can choose a finite number of points $x_1,\ldots, x_N$ such that 
$$
\partial S\subset\bigcup_{i=1}^N \Omega^{x_i}\doteq {\cal D}\,.
$$ 
By denoting with $\alpha_i$ the constant $\alpha$ from Step~1 corresponding to the point $x_i\in\partial S$, for each $i\in\{1,\ldots,N\}$, we re--label the sets as follows
$$
\Omega_{\alpha_{i},i}\doteq \Omega^{x_i}\,,
$$
and we order the collection $\{\Omega_{\alpha_{i},i}\}$ with lexicographic order, i.e.
$$
(\alpha_{i},i)\prec(\alpha_{j},j)\quad\mbox{ iff }\quad 
\mbox{ either }~\alpha_{i} < \alpha_{j} ~\mbox{ or } ~\alpha_{i} = \alpha_{j} \mbox{ and } i<j\,.
$$
The advantage of this choice is that we can easily prove that points of $\partial^- \Omega^{x_i}$, i.e. of the ``lower'' boundary of $\Omega^{x_i}$, belong either to $\partial {\cal D}$ or to another domain $\Omega_{\alpha_{j},j}$ with larger index. Namely, fix $i\in\{1,\ldots,N\}$ and let $z\in\partial^- \Omega^{x_i}$. If $z\notin\partial {\cal D}$, then there holds
$$
z\in\Omega^{x_{j^*}}\,,\qquad\qquad j^*\doteq\max\big\{j~;~z\in\Omega^{x_j}\big\}\,, \qquad\qquad j^*\neq i\,.
$$
Hence, if we assume that
$
(\alpha_{i},i)\not\prec(\alpha^*,{j^*})
$, 
where we have set for simplicity $\alpha^*\doteq\alpha_{{j^*}}$, we would have in particular that $\alpha_{i}\geq \alpha^*$ and therefore
$$
d(z,\R^d\setminus S)=\alpha_{i}\geq \alpha^*
\qquad\Longrightarrow\qquad
z\in S_{\alpha^*}\,.
$$
This means $z\notin \Omega^{x_{j^*}}=\Gamma^{x_{j^*}}\setminus  S_{\alpha^*}$, i.e., a contradiction. Thus, it must be $(\alpha_{i},i)\prec(\alpha^*,{j^*})$.

\medskip

We now claim that defining a feedback control $U$ by setting
$$
U(\xi)\doteq u_{x_i} \qquad\qquad\mbox{for all }~~\xi\in\Omega_{\alpha_{i},i}\setminus \bigcup_{(\alpha_{i},i)\prec (\alpha_{j},j)}\Omega_{\alpha_{j},j}
$$
we obtain a patchy feedback on ${\cal D}$ with the required properties. 

\n Indeed, if $z \in\partial^+\Omega^{x_i}$ for a certain $i=1,\ldots,N$ and if we denote with $\tilde\ve_i$ the constant $\tilde\ve$ from Step~1 corresponding to the point $x_i\in\partial S$, then we have $|f(z,u_{x_i})-f(x_i,u_{x_i})|<\tilde\ve_i/4$ by construction. 
Observing that for any pair of vectors $v,v'\in\R^d$ and any radius $\ell$, one always has
$$
|v-v'|<\,{\ell\over 2}\qquad\Longrightarrow\qquad
{\cal W}(v',\ell/2)\subseteq{\cal W}(v,\ell)\,,
$$
we obtain that  $z+cf(z,u_{x_i})\in z+{\cal W}(f(x_i,u_{x_i}),\tilde\ve_i/2)$ for some $0<c\leq1$, proving the inward pointing condition at $z$. This proves that $U$ is a patchy feedback.

%
%
\n Next, we need to prove that there exists $\tilde r\in\,]0,r_o[$ such that $Q(S,\tilde r)\subseteq{\cal D}$. If this was not the case, then we would have a sequence of points $\xi_n$ on $\partial{\cal D}\cap S$ such that $d(\xi_n,\overline{\R^d\setminus S})\to 0$. 
By possibly taking a subsequence, which we do not relabel, $\xi_n\to\bar \xi\in\partial S$. But since $\partial {\cal D}$ is closed, it must be $\bar \xi\in\partial S\cap\partial {\cal D}$, and then we would have a contradiction, because it would mean
$$
\bar \xi~\in~ \partial S\setminus\left(\bigcup_{i=1}^N \Omega^{x_i}\right)\neq\emptyset\,.
$$

\n Finally, inclusion $\Omega^{x_i}\subseteq x_i+\rho_{x_i} B_d$ guarantees that~\eqref{eq:ipc_constraint2} holds in $\Omega^{x_i}$, for each $i\in\{1,\ldots,N\}$, with a positive constant $\mu_{x_i}$. Hence,~\eqref{eq:inward_bdry} holds by simply taking $\mu\doteq\min\{\mu_{x_1}/2,\ldots,\mu_{x_N}/2\}$. Since we had already observed that $\rho_x$ in~\eqref{eq:inclusion_lemma2} can be chosen so that $\diam\,\Omega^x\leq\lambda$, for any fixed constant $\lambda>0$, the proof is complete.~~$\diamond$

\section{Proof of the main result}

Let $\delta>0$ be fixed. We want to construct a patchy feedback control $U\colon{\cal D}\mapsto{\bf U}$, defined on a domain ${\cal D} \supseteq S\setminus\Sigma^\delta$, such that every Carath\'eodory solution to~\eqref{eq:system_cl} with $x(0)\in S\setminus\Sigma^\delta$ are steered to $\Sigma^\delta$ in finite time. Let $\gamma>0$ be a constant such that $4\gamma\leq\delta$.

\medskip

\n {\bf Step~1}. We start with the construction of the patchy feedback close to the boundary $\partial S$ of the constraint. By applying Lemma~\ref{lem:patchy_from_wedged}, there exist a constant $\tilde r\in\,]0,r_o[$, with $r_o>0$ the value found in Lemma~\ref{lem:CRS1}, a constant $\mu>0$ and a patchy feedback control $U_o=\big(U_o,(Q^o_\alpha,q^o_\alpha)_{\alpha\in\A_o}\big)$, defined on a domain ${\cal D}_o\supset Q(S, \tilde r)$, such that~\eqref{eq:inward_bdry} holds for all $x\in Q(S, \tilde r)$. Since $Q(S,\tilde r)$ is compact, it is not restrictive to assume that the set of indices $\A_o$ is finite, i.e., $\A_o=\{1,\ldots,N_o\}$ for a suitable $N_o\in\N$.

\n Notice that~\eqref{eq:inward_bdry} is equivalent to say that there exists $\ve>0$ such that
$$
y+{\cal W}(f(x,U(x)),\ve)\subset S_{r(x)}\qquad\qquad~~~ \forall~y\in \{x+\ve B_d\}\cap S_{r(x)}\,,~~\forall~x\in Q(S, \tilde r)\,.
$$
Bearing~\eqref{eq:cones_to_wedged} in mind, this also means that
$$
y+{\cal W}(-f(x,U(x)),\ve)\subset \overline{\R^d\setminus S_{r(x)}}\qquad~~~ \forall~y\in \{x+\ve B_d\}\setminus\overset{\,\circ}{S}_{r(x)}\,,~~\forall~x\in Q(S, \tilde r)\,.
$$
Hence, for every $x\in Q(S, \tilde r)$, by applying Lemma~\ref{lem:CRS4} to the vector $v=f(x,U(x))$, to the wedged sets $S_{r(x)}$ and to the point $x\in\partial S_{r(x)}$, we deduce the existence of a neighborhood ${\cal N}_x$ of $x$ such that
\begin{equation}\label{eq:dist_decr}
\nabla\Delta_{S_{r(x)}}(y)\cdot f(x,U(x)) \leq -\ve<0\,,
\end{equation}
whenever $y\in{\cal N}_x$ is a point of differentiability of the map $\xi\mapsto\Delta_{S_{r(x)}}(\xi)$.

\medskip

\n {\bf Step~2}. Next, observe that, fixed any $\xi\in S_{\tilde r} \setminus \Sigma$, the assumption of open loop $S$--constrained controllability to $\Sigma$ ensures the following. There exist a measurable open loop control $u_{\xi}(\cdot)$ and a time $T_\xi\geq 0$ such that the Carath\'eodory solution $x(\cdot)$ to~\eqref{eq:system}, starting from $\xi$ with control $u_{\xi}$, satisfies
$$
x(t)\in S\qquad\forall~t\in[0,T_\xi]\,,\qquad\qquad\qquad x(T_\xi)\in\Sigma\,.
$$
By possibly reducing $T_\xi$ and redefining $u_{\xi}(\cdot)$ on some subinterval of $[0,T_\xi]$, it is not restrictive to assume that $x(t')\neq x(t)$ whenever $t'\neq t$ in $[0,T_\xi]$.

Owing to Corollary~\ref{cor:stick_inside} with $r=\rho_{\xi}\doteq\min\left\{\,{\gamma\over 2C_{T_\xi}},\tilde r\right\}$, where $C_{T_\xi}$ is the constant in the corollary corresponding to the time $T_\xi$, we can find another control $\bar u_{\xi}(\cdot)$ defined on $[0,T_\xi]$ so that the corresponding solution $\bar x(\cdot)$ emanating from $\xi$ satisfies
\begin{equation}\label{eq:approx}
\bar x(t)\in S_{\rho_{\xi}}\qquad\forall~t\in[0,T_\xi]\,,\qquad\qquad\qquad \bar x(T_\xi)\in\Sigma^{\gamma/2}\,.
\end{equation}
Since the set of piecewise constant admissible controls is dense in the set of all controls, we can also find a piecewise constant control $\tilde u_\xi$ which approximates $\bar u_\xi$ so that the corresponding solution $\tilde x(\cdot)$ starting from $\xi$ satisfies
\begin{equation}\label{eq:approx2}
\tilde x(t)\in S_{\rho_{\xi}/2}\qquad\forall~t\in[0,T_\xi]\,,\qquad\qquad\qquad \tilde x(T_\xi)\in\Sigma^\gamma\,.
\end{equation}
Then, we are in the position to apply Lemma~\ref{lem:patchy_tubes} with $\ve=\min\{\rho_{\xi}/4,\gamma\}$ to obtain an open domain $\Gamma_{\xi}$ and a patchy feedback control $U_\xi=\big(U_\xi,(Q^\xi_\alpha,q^\xi_\alpha)_{\alpha\in\A_\xi}\big)$, defined on a domain ${\cal D_\xi}\supseteq\Gamma_\xi$, such that properties {\it (i)}--{\it (iv)} of the lemma are satisfied. In particular, there holds
$$
\Gamma_{\xi}\subseteq S_{\rho_{\xi}/4}\subset S\,,
$$
and every Carath\'eodory solution to~\eqref{eq:system_cl}, with initial datum in $\Gamma_\xi$ and control $U_\xi$, remains inside $\Gamma_\xi\subseteq S$ and eventually reaches $\tilde x(T_\xi)+\ve\overline{B_d}\subseteq\Sigma^{2\gamma}$. Moreover, by exploiting the inward--pointing properties {\it (ii)} across $\partial\Gamma_\xi$ to possibly replace some domains $Q^\xi_\alpha$ with $Q^\xi_\alpha\cap\Gamma_\xi$ and by slightly modifying the patches which intersect $\Sigma^{2\gamma}$, 
%
%
we can assume that
\begin{equation}\label{eq:strict_tubes}
\bigcup_{\alpha\in\A_\xi} Q_\alpha^\xi \setminus \Sigma^{3\gamma}=\Gamma_\xi\setminus \Sigma^{3\gamma}
\subseteq S_{\rho_{\xi}/4}\,.
\end{equation}

\medskip

\n {\bf Step~3}. In the previous step, we have constructed for every $\xi\in S_{\tilde r} \setminus \Sigma$ a patchy feedback control $U_\xi$, defined on an open domain $\Gamma_\xi$, so that properties {\it (i)}--{\it (iv)} of  Lemma~\ref{lem:patchy_tubes} are satisfied and~\eqref{eq:strict_tubes} holds. Then, by compactness of the set $S_{\tilde r} \setminus \Sigma^{3\gamma}$, we can find a finite number of points $\xi_1,\ldots,\xi_M$ such that 
$$
S_{\tilde r} \setminus \Sigma^{3\gamma}\subseteq\bigcup_{j=1}^M \Gamma_{\xi_j}\,.
$$
Since for $j=1,\ldots,M$ each $\overline{\Gamma}_{\xi_j}$ is compact, it is not restrictive to assume that the corresponding patchy feedbacks $U_{\xi_j}=\big(U_{\xi_j},(Q^j_\alpha,q^j_\alpha)_{\alpha\in\A_j}\big)$ have finite sets of indices $\A_j=\{1,\ldots,N_j\}$.
By defining a new set of indices ${\cal A}\doteq\{1,\ldots,N_o,N_o+1,\ldots,N_o+N_1,\ldots,\sum_{j=0}^M N_j\}$, by setting
$$
\forall~\alpha\in\A~\mbox{ s.t. }~\sum_{j=0}^{k-1}N_j<\alpha\leq \sum_{j=0}^kN_j\,,
\qquad
Q_\alpha\doteq Q^k_{\alpha-\sum_{j=0}^{k-1}N_j}\,,
\qquad
q_\alpha \doteq q^k_{\alpha-\sum_{j=0}^{k-1}N_j}\,,
$$
%
and
$$
U(x) \doteq  q_\alpha\qquad \hbox{if}\qquad x \in
Q_\alpha \setminus \bigcup_{\beta > \alpha} Q_{\beta}\,, \qquad  \alpha\in\A\,,
$$
it is easy to verify that $U=\big(U,(Q_\alpha,q_\alpha)_{\alpha\in\A}\big)$ is a patchy feedback control on the domain
$$
{\cal D}\doteq \left({\cal D}_o\cup\bigcup_{j=1}^M\Gamma_{\xi_j}\right) \setminus \Sigma^{3\gamma}\supset
S \setminus \Sigma^{\delta}\,.
$$
Indeed, the inward--pointing condition to be verified by $f(z,U(z))$ at points $z$ of $\partial Q_\alpha\setminus\bigcup_{\beta>\alpha}  Q_{\beta}$, $\alpha\in\A$, is part of the corresponding conditions verified by the patchy feedbacks $U_o,U_{\xi_1},\ldots,U_{\xi_M}$.

We claim that such patchy feedback has the required properties, i.e. that trajectories of~\eqref{eq:system_cl} corresponding to the control $U$ and starting from a point in $S\setminus \Sigma^\delta$ do not exit from $S$ and eventually enter $\Sigma^\delta$. 
The proof of this claim is given in the next two steps and completes the proof of Theorem~\ref{thm:practical_stab_constr}.

\medskip

\n {\bf Step~4}. We prove that any trajectory $x(\cdot)$ of~\eqref{eq:system_cl}, corresponding to the control $U$ and such that $x(0)\in S\setminus\Sigma^\delta$, remains inside $S$ for all times $t\geq 0$ in its maximal domain of existence $[0,T_{max}[$. 

\smallskip

\n First, observe that it is enough to prove the property for trajectories $x(\cdot)$ with $x(0)$ in the interior of $S$. Indeed, if $x(0)\in \partial S$, then there exists a small $\tau>0$ such that in $]0,\tau]$ a solution exists, because the vector field is patchy, and belongs to the interior of $S$, because of~\eqref{eq:inward_bdry}. Hence, by applying  the result on the interior of $S$ to the solution of~\eqref{eq:system} with initial datum $x(\tau)$, one concludes that the property holds also for trajectories starting from the boundary $\partial S$.

\smallskip

\n Now, take $x(0)\in \overset{\,\,\circ}{S}$ and assume by contradiction there exists $\bar t\in\,]0,T_{max}[$ such that $x(\bar t)\in\partial S$ but $x(s)$ belongs to the interior of $S$ for $0\leq s<\bar t$. Then,
$$
x(\bar t)\in Q_{\bar \alpha}\,,
\qquad\mbox{where}\qquad
\bar \alpha\doteq\max\{\alpha\in\A~;~x(\bar t)\in Q_\alpha\}\,.
$$
Actually, recalling the construction of the domains $\Gamma_{\xi_1},\ldots,\Gamma_{\xi_M}$ in steps~2--3, one has that
%
%
$$
\bigcup_{j=1}^M\Gamma_{\xi_j}\subseteq S_{\bar \rho/4}\,,\qquad\qquad
\bar\rho\doteq\min\{\rho_{\xi_1},\ldots,\rho_{\xi_M}\}\,.
$$
Hence, $Q_{\bar \alpha}$ must be one of the domains used in $U_o$, i.e. $\bar \alpha\in\A_o$, which in turn implies that $f(\cdot,U(x(\bar t)))=f(\cdot,q_{\bar \alpha})$ satisfies~\eqref{eq:dist_decr}. More precisely, observe that, since $Q_{\bar \alpha}$ is an open set and $U$ is a patchy control,  there exists $\tau\in [0,\bar t[$ such that for $s\in \,]\tau,\bar t]$ one has 
$$
x(s)\in Q_{\bar \alpha}\,,\qquad\qquad \dot x(s)=f(x(s),q_{\bar \alpha})\,.
$$
Now assume that $\Delta_S$ is differentiable in $x(s)$ for a.e. $s\in \,]\tau,\bar t]$. Then, by combining the local decrease property~\eqref{eq:dist_decr}, applied in a neighborhood of $x(\bar t)\in\partial S$, with the boundedness of $\nabla \Delta_S$ and the continuity of $f$ in $Q_{\bar \alpha}$, there exists $\tau'\in\,[\tau,\bar t[$ such that for a.e. $s\in \,[\tau',\bar t]$ there holds
\begin{align*}
\nabla\Delta_S(x(s))\cdot f(x(s),q_{\bar \alpha})&\leq\nabla\Delta_S(x(s))\cdot f(x(\bar t),q_{\bar \alpha})\\
&~~~~~~~~~~~~~~+\Big|\nabla\Delta_S(x(s))\Big|~\Big| f(x(\bar t),q_{\bar \alpha})-f(x(s),q_{\bar \alpha})\Big| \leq -\ve/2<0\,.
\end{align*}
In turn, this implies
\begin{align}\label{eq:contrad_ineq}
d(x(\bar t),\overline{\R^d\setminus S})-d(x(\tau'),\overline{\R^d\setminus S})&=-\Delta_S(x(\bar t))+\Delta_S(x(\tau'))\nonumber\\
&= -\int_{\tau'}^{\bar t}\nabla\,\Delta_S(x(\sigma))\cdot \dot x(\sigma)\,d\sigma\geq \,{\ve\over 2}\,(\bar t-\tau')>0\,,
\end{align}
which yields a contradiction since $d(x(\bar t),\overline{\R^d\setminus S})=0$ and $x(\tau')\in\overset{\,\,\circ}{S}$.
Notice that we can reach the same conclusion even if $\Delta_S$ is not differentiable along the trajectory $x(\cdot)$ on a set ${\cal I}\subseteq\,]\tau,\bar t]$ of positive measure. Indeed, denote with $\tau'\in\,]\tau,\bar t]$ a value such that $x(s)\in{\cal N}_{x(\bar t)}$ for $s\in\,[\tau',\bar t]$, where ${\cal N}_{x(\bar t)}$ is the neighborhood of $x(\bar t)$ in which~\eqref{eq:dist_decr} holds. Then, it is always possible to find an arbitrarily close curve $x_\alpha(\cdot) \doteq x(\cdot)+\alpha$, $\alpha\in\R^d$ with $|\alpha|\ll 1$, where $\Delta_S$ is differentiable a.e., since otherwise the Lipschitz continuous function $\Delta_S$ would be not differentiable on a subset of positive measure of the following neighborhood of $x(\cdot)_{|{[\tau',\bar t]}}$
$$
\big\{ x(s)+ \rho B_d~;~ s\in\,[\tau',\bar t],\,\rho\ll 1\big\}\,.
$$
Hence, we can repeat the computation above for the variation of $\Delta_S$ along such a curve $x_\alpha(\cdot)$: by $\dot x_\alpha = \dot x$ and by the arbitrary closeness of $x_\alpha$ to $x(\cdot)$, we thus reach a contradiction also in this case.

\medskip

\n {\bf Step~5}. Finally, it remains to prove that all trajectories are steered to $\Sigma^\delta$. Fix an initial datum $x(0)\in S\setminus\Sigma^\delta$ and denote with $x(\cdot)$ the trajectory of~\eqref{eq:system_cl} corresponding to the patchy feedback $U$. If $x(0)\in S_{\tilde r}\setminus\Sigma^\delta$, then $x(0) \in\Gamma_{\xi_j}$ for some $j\in\{1,\ldots, M\}$. By property~{\it (iii)} in Lemma~\ref{lem:patchy_tubes}, each domain $\Gamma_{\xi_j}$ is positively invariant for the dynamics with control $U_{\xi_j}$, hence there exist $n\geq 1$, times $t_o=0<t_1<\ldots<t_n$ and integers $\ell_1<\ldots<\ell_n<\ell_{n+1}$ in $\{1,\ldots, M\}$ such that 
$$
x(t)\in\Gamma_{\xi_{\ell_i}}\qquad\qquad \forall~t\in(t_{i-1},t_{i}]\,,~~~i=1,\ldots,n
$$
and 
$$
x(t)\in\Gamma_{\xi_{\ell_{n+1}}}\qquad\qquad \forall~t>t_n\,.
$$
%
%
Observing that, for $t>t_n$, the trajectory $x(\cdot)$ is a Carath\'eodory solution to the Cauchy problem
$$
\dot y = f(y,U_{\xi_{\ell_{n+1}}}(y))\,,\qquad\qquad y(t_n)=x(t_n)\,,
$$
we can use property~{\it (iv)} of Lemma~\ref{lem:patchy_tubes} to compare $x(\cdot)$ with the trajectory $\bar x(\cdot)$ emanating from $\xi_{\ell_{n+1}}$, which by construction is defined on some interval $[0,\overline{T}]$ and satisfies $\bar x(\overline{T})\in\Sigma^\gamma$. Namely, we can conclude that the trajectory $x(\cdot)$ exists until a certain time $t'>t_n$ such that $x(t')\in \bar x(\overline{T})+\ve B_d\subseteq\Sigma^{2\gamma}\subset\Sigma^\delta$.

%
%

Assume now that the initial datum $x(0)$ is in $Q(S,\tilde r)\setminus\Sigma^\delta$ and denote again the corresponding solution with $x(\cdot)$ and its maximal domain of existence with $[0,T_{max}[$. Since we have proved in Step~4 that $x(\cdot)$ cannot leave $S$, and since $\partial{\cal D}\cap S\subseteq\Sigma^\delta$, either the trajectory eventually reaches $\Sigma^\delta$ or $T_{max}=+\infty$ and $x(\cdot)$ remains in $S\setminus\Sigma^\delta$ for all times. Having already proved that $x(\tau)\in\bigcup_{j=1}^M\Gamma_{\xi_j}$ for some $\tau\geq 0$ implies stabilization to $\Sigma^\delta$, we have only to exclude the case 
$$
x(t)\in Q(S,\tilde r)\setminus\left(\Sigma^\delta~\cup~\bigcup_{j=1}^M\Gamma_{\xi_j}\right)\qquad\qquad\forall~t\geq 0\,.
$$ 
But this can be ruled out by repeating the arguments used in Step~4. Fixed any $r\in\,]0,\tilde r[$, property~\eqref{eq:inward_bdry} ensures that trajectories starting from a point in $\partial S_r$ will belong to $\overset{\,\,\circ}{S}_r$ for small positive times. Moreover, fixed any interval $[\tau,\tau'[\,\subset[0,+\infty[$, if $x(s)\in\overset{\,\,\circ}{S}_r$ for $s\in[\tau,\tau'[$, then it cannot be $x(\tau')\in\partial S_r$. Indeed, if that were the case, we would recover a contradiction with the local decrease property~\eqref{eq:dist_decr} of the function $z\mapsto\nabla \Delta_{S_r}(z)$, in a neighborhood of $x(\tau')$. In conclusion, $d(x(t),\overline{\R^d\setminus S})$ must increase with time and the trajectory eventually reaches either $\Sigma^\delta$ or $S_{\tilde r}\subseteq\bigcup_{j=1}^M\Gamma_{\xi_j}$.~~$\diamond$

\section{Extensions and remarks}\label{sec:final}

\subsection{An explicit example}
Let us consider the following control system in $\R^2$ (a perturbed harmonic oscillator from~\cite{CS1})
\begin{equation}\label{eq:explicit_eq}
\dot x = \left(
\begin{array}{cc}
0 & 1\\
-1 & 0
\end{array}
\right)\, x + u\, x\,,
\qquad\qquad
u\in [-1,1]\,.
\end{equation}
As constraint and target sets, let us consider respectively
$$
S=\{x\in\R^2~;~ 1\leq |x|\leq 2\}\,,
\qquad\qquad
\Sigma = \{(2,0)\}\,.
$$
It is easy to see that such system satisfies all the assumptions of Theorem~\ref{thm:practical_stab_constr}, including the open--loop $S$ constrained controllability to $\Sigma$. Let $\delta\in\,]0,1[$ be fixed and let us construct a patchy feedback which steers every trajectory into the neighborhood $\Sigma^\delta=B((2,0),\delta)$.

\begin{figure}
\centering
\psfrag{t}{$\scriptstyle \bar\theta$}
\psfrag{O1}{$\Omega_1$}
\psfrag{O}{$\Omega_0$}
\psfrag{S}{$S'$}
\psfrag{B}{$\Sigma^\delta$}
\includegraphics[width=0.35\textwidth]{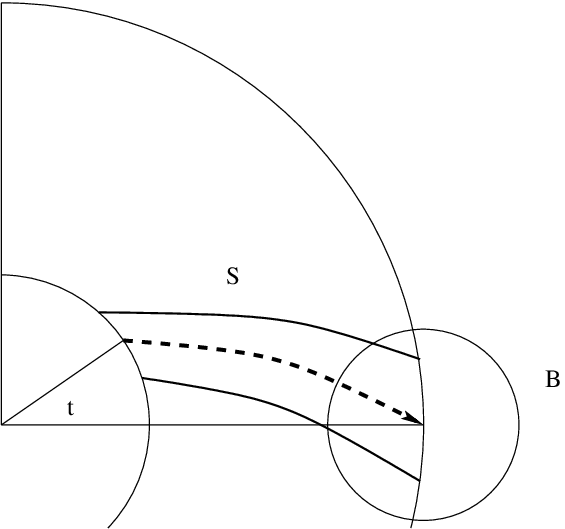}
\qquad \qquad \qquad
\includegraphics[width=0.25\textwidth]{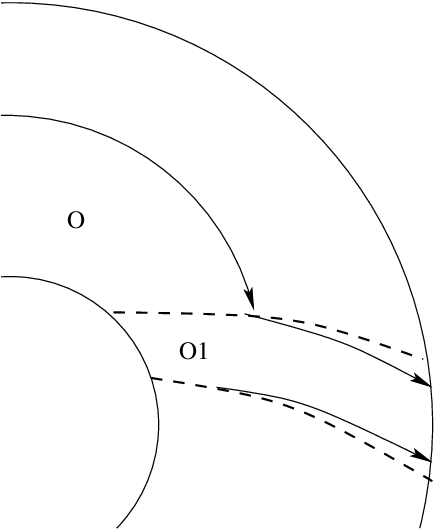}
\caption{(Left) Tube--like domain $S'$ around the dashed trajectory $\bar x(\cdot)$. (Right) Some trajectories steered to $\Sigma^\delta$ by the patchy control $U$.}
\label{fig:ex1}       
\end{figure}

Observing that we can rewrite~\eqref{eq:explicit_eq} in polar coordinates as
$$
\dot \rho = u\,\rho\,,
\qquad\qquad
\dot \theta = -1\,,
$$
one can check that the trajectory $\bar x(\cdot)$ corresponding to $u\equiv 1$, and starting from the point $(\bar\rho,\bar\theta)=(1, \ln 2)$, arrives exactly at $\Sigma$ in time $T=\ln 2$. Now, by defining (cf. Figure~\ref{fig:ex1}--left)
$$
S'\doteq \left\{(\rho,\theta)\in S~;~\rho\in\,]1,2[\,,~\theta-\bar\theta\in\, \left]-\theta_0-\left(1+\,{\theta_0\over ln 2}\right)\ln\rho,\theta_0-\left(1-\,{\theta_0\over ln 2}\right)\ln\rho \right[~
\right\}
$$
for a suitable fixed $\theta_0$ small enough so to have 
$\{(\rho,\theta)\in S~;~\rho=2\,,\theta\in [-2\theta_0,2\theta_0]\}\subseteq \Sigma^\delta$ (for instance $\theta_0\approx {1\over 4}\arcsin{\delta\over 4}$ has this property), we obtain a tube--like open set $S'$ around the trajectory $\bar x(\cdot)$ such that the set is positively invariant for the dynamics with constant control $u\equiv 1$ and any trajectory starting inside $S'$ arrives into $\Sigma^\delta$ in time $T<\ln 2$.
Hence, we can set as in Figure~\ref{fig:ex1}--right
$$
\Omega_0 \doteq \R^2\qquad U_0\equiv 0\,,
\qquad\qquad\qquad
\Omega_1 \doteq S'\qquad U_1\equiv 1\,,
$$
to obtain a $S$--constrained patchy control $(U, (\Omega_\alpha,U_\alpha)_{\alpha=0,1})$ which steers every trajectory to $\Sigma^\delta$.

Observe that, differently from the patchy feedback constructed in Theorem~\ref{thm:practical_stab_constr}, $U$ does not give a strictly inward pointing vector field at points of $\partial S$ and thus $U$ is less robust w.r.t. external disturbances: trajectories which travels along the boundary might violate the constraint whenever a very small perturbation pointing outside $S$ occurs. To recover robustness close to $\partial S$, we must follow more closely the proof of Theorem~\ref{thm:practical_stab_constr}. 

\begin{figure}
\centering
\psfrag{O1}{$\!\!\!\!\!\!\!\!\Omega_1$}
\psfrag{O}{$\!\!\!\!\!\!\!\Omega_0$}
\psfrag{O2}{$\!\!\!\!\!\!\!\Omega_2$}
\psfrag{O3}{$\Omega_3$}
\psfrag{S}{$\!\!\!\!\!S'$}
\psfrag{S2}{$S''$}
\psfrag{B}{$\!\!\!\Sigma^\delta$}
\includegraphics[width=0.45\textwidth]{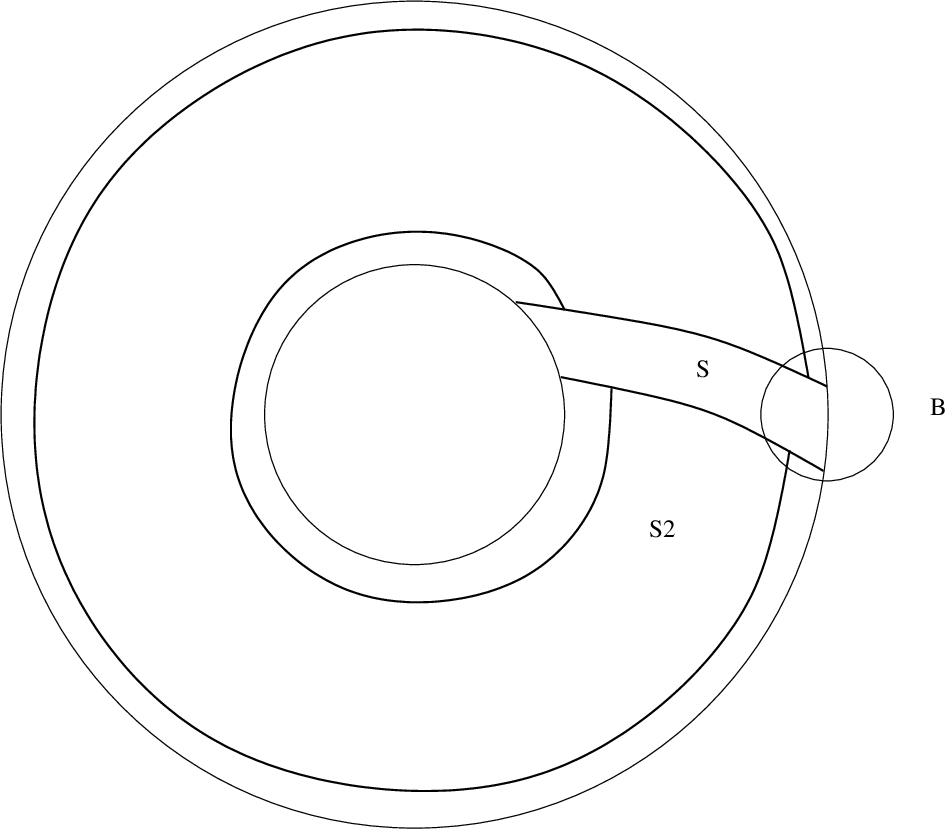}
\qquad
\includegraphics[width=0.4\textwidth]{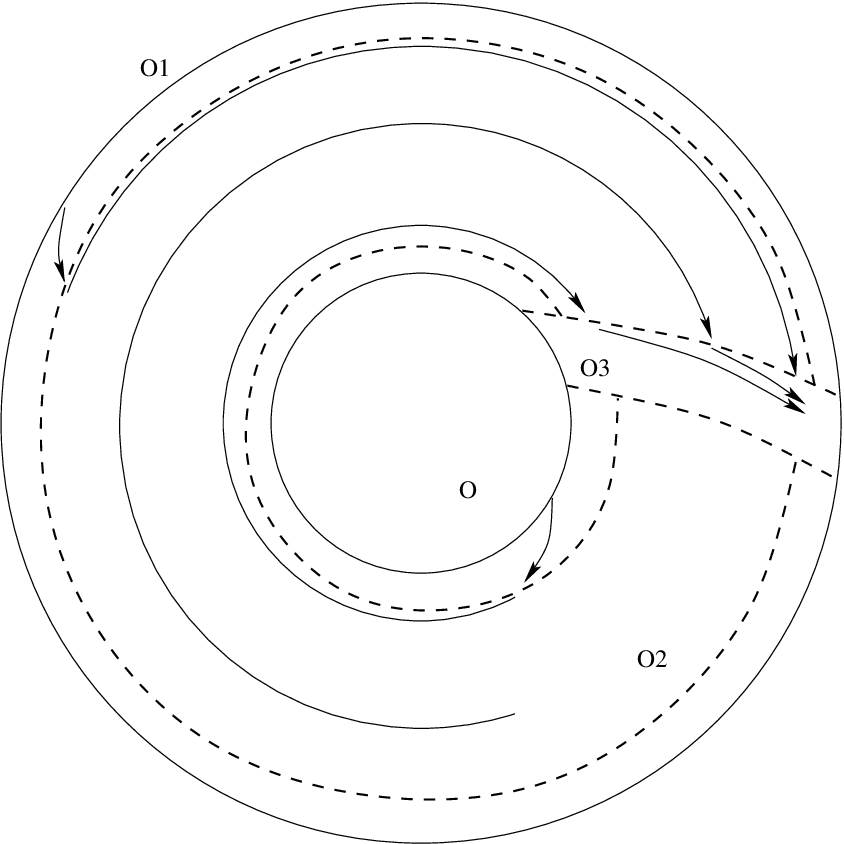}
\caption{(Left) Tube--like domain $S''$ around circular trajectories with $u\equiv 0$ (of increasing size as $\theta$ runs counterclockwise). (Right) Some trajectories steered to $\Sigma^\delta$ by the patchy control $V$.}
\label{fig:ex2}       
\end{figure}

\n Fixed a small parameter 
$\tau\in\,]0,e^{\theta_0}-1[\,$, 
we define the following set (see Figure~\ref{fig:ex2}--left)
$$
S''\doteq\left\{(\rho,\theta)\in S~;~\theta\in\,]0,2\pi[\,,~\rho\in\,\left]1+{\tau\over 2}\Big(1+\Big\{{\theta-\bar\theta\over 2\pi}\Big\}\,\Big),2-{\tau\over 2}\Big(1+{\theta\over 2\pi}\Big)\right[~\right\}\,,
$$
where $\{x\}$ is the fractional part of the real number $x$ (i.e., $\{x\}=x-\lfloor x\rfloor$). \footnote{Observe that the boundary curves $\rho_2(\theta)=2-{\tau\over 2}\Big(1+{\theta\over 2\pi}\Big)$ and $\rho_1(\theta)=1+{\tau\over 2}\Big(1+\Big\{{\theta-\bar\theta\over 2\pi}\Big\}\,\Big)$ of $S''$ are simply two curves which approach the boundaries $\rho=2$ and $\rho=1$ of $S$, respectively, as $\theta$ varies counterclockwise. The curve $\rho_1$ is additionally shifted of an angle $\bar\theta$ to ensure that its endpoints belongs to $S'$.}
It is easy to verify that $S''$ is positively invariant for the dynamics with constant control $u\equiv 0$ and that trajectories starting inside $S''\setminus S'$ arrive either into $\Sigma^\delta$ or into $S'$ in time $T<2\pi$. 
Therefore, by setting as in Figure~\ref{fig:ex2}--right
$$
\begin{array}{lll}
\Omega_0 \doteq \{x\in\R^2~;~ 0<|x|<1+\tau\} &\qquad\qquad& V_0\equiv 1\,,\\&&\\
\Omega_1 \doteq \{x\in\R^2~;~ |x|>2-\tau\} &\qquad\qquad& V_1\equiv -1\,,
\end{array}
$$
$$
\Omega_2 \doteq S''\qquad V_2\equiv 0\,,
\qquad\qquad\qquad
\Omega_3 \doteq S'\qquad V_3\equiv 1\,,
$$
we obtain a $S$--constrained patchy feedback $(V, (\Omega_\alpha,V_\alpha)_{\alpha=0,1,2,3})$ which steers every trajectory to $\Sigma^\delta$ and $V$ is also robust in a neighborhood of $\partial S$ w.r.t. to any small internal or external perturbation of the dynamics, because it uses there one of the control values satisfying {\bf (S2)}.

\subsection{$S$--restricted dynamics}

Our main result has been stated and proved for a control system~\eqref{eq:system} whose vector field $f$ is defined in the whole $\R^d\times{\bf U}$. In many applications to economy and engineering, though, the dynamics could have no meaning or even break down completely when $x\notin S$. As such we want to stress that Theorem~\ref{thm:practical_stab_constr} can be also applied to the case of a control system~\eqref{eq:system} whose dynamics is given by a function $f\colon S\times {\bf U}\to\R^d$ not defined for $x\notin S$.

\smallskip

Indeed, in the same spirit of~\cite{CRS, CS1}, we can extend $f$ to a globally Lipschitz continuous function $\tilde f\colon \R^d\times {\bf U}\to\R^d$ by defining $\tilde f=(\tilde f_1,\ldots,\tilde f_d)$ as follows
\begin{equation}\label{eq:extended_vf}
\tilde f_i(x,u)\doteq \min_{y\in S}\big\{ f_i(y,u)+L_f|x-y|\big\}\,.
\end{equation}
Now, assuming that the vector field $f$ satisfies {\bf (F1)}, {\bf (F2)} and {\bf (F3)} on its domain $S$, then also $\tilde f$ satisfies {\bf (F1)} and {\bf (F2)} on the whole $\R^d$. However, $\tilde f$ might fail to satisfy {\bf (F3)} outside $S$. Luckily, this is not a problem for our result: we needed convexity of $f(x,{\bf U})$ only to apply Lemma~\ref{lem:track} (tracking lemma) close to the boundary $\partial S$. As such, {\bf (F3)} is only necessary inside $S$, and the rest of our construction can be applied to the dynamics 
\begin{equation}\label{eq:system_ext}
\dot x=\tilde f(x,u)\,,
\end{equation}
with no significant change. 


\subsection{Unbounded constraints}

Another assumption that can be relaxed in our main result is the compactness of the constraint set $S$. Indeed, as it has been done for sample--and--hold trajectories in~\cite{CS1}, one can require $S$ to only be a closed set. In the latter case, it is possible to prove the following

\begin{theorem}\label{thm:unbdd_stab_constr} Let $S$ be a set satisfying {\bf (S1)} and {\bf (S2)}, except for compactness in {\bf (S1)} replaced by closedness, and let $\Sigma$ be any closed set such that $S\cap \Sigma\neq\emptyset$. Assume that for all bounded sets of initial data ${\cal B}$, the trajectories of the system~\eqref{eq:system} starting from ${\cal B}$ are open loop controllable to $\Sigma$ remaining inside $S$. Then, for every fixed bounded set ${\cal B}$ there exists a patchy feedback control $U=U_{\cal B}(x)$ which makes~\eqref{eq:system} practically stable to $\Sigma$ subject to the constraint $S$.
\end{theorem}

To prove the theorem above, notice that the results about wedged sets given in Section~\ref{sec:prelim} still hold under the relaxed assumptions of Theorem~\ref{thm:unbdd_stab_constr} (see~\cite{CS1}). Therefore, for the fixed bounded set ${\cal B}$, one can simply take a large enough ball $K=k B_d$ such that ${\cal B}\subseteq K$ and $(\overline{S\cap K})\cap\Sigma\neq\emptyset$, 
and apply Theorem~\ref{thm:practical_stab_constr} to the smaller constraint set $S'\doteq\overline{S\cap K}$. In this way, for every $\delta>0$ we obtain a patchy feedback control such that trajectories starting from ${\cal B}$ remain inside $S$ for all positive times and eventually reach $(S\cap K)\cap\Sigma^\delta\subseteq\Sigma^\delta$, as required.

\subsection{Robustness}\label{sec:rob}

One of the main advantages of using patchy controls $U(x)$ and Carath\'eodory solutions for~\eqref{eq:system} over allowing arbitrarily discontinuous controls and weaker concepts of solutions like sample--and--hold trajectories (see~\cite{CRS,CS1,CS2}), is that stronger robustness properties can be proved with almost no efforts.

\n It is indeed well known that whenever the vector field
$$
g(x)= f(x,U(x))\,,
$$
is a patchy vector field in the sense of Definition~\ref{defn:patchy_vf}, then the set of Carath\'eodory solution is robust w.r.t. both internal and external perturbations (see~\cite{AB1,AB2,ABsurvey}) without any additional assumption on the feedback control $U$. This represents a noticeable improvement compared to the construction obtained in~\cite{CRS, CS1, CS2} through sample--and--hold and Euler solution, which only achieves the same robustness by requiring a ``reasonable uniformity'' in the time discretization.

\smallskip

This additional robustness holds also in the case of a constrained dynamics. Namely, the following theorem holds. 

\begin{theorem}\label{thm:stab_constr_robust} Assume that the system~\eqref{eq:system} satisfies open loop $S$--constrained controllability to $\Sigma$, where $S$ is a set satisfying {\bf (S1)} and {\bf (S2)} and $\Sigma$ is any closed set such that $S\cap \Sigma\neq\emptyset$. Then, for every $\delta>0$ there exist $T>0$, $\chi>0$ and a patchy feedback control $U=U(x)$, defined on an open domain ${\cal D}$ with $S\setminus\Sigma^\delta\subseteq {\cal D}$, so that the following holds. Given any pair of maps $\zeta\in\BV([0,T], \R^d)$ and $d\in\Lsp^1([0,T], \R^d)$ such that
$$
\|\zeta\|_{\BV}\doteq\|\zeta\|_{\Lsp^1}+\mathrm{Tot. Var.}(\zeta)<\chi\,,\qquad\qquad
\|d\|_{\Lsp^1}<\chi\,,
$$
and any initial datum $x_o\in S\setminus\Sigma^\delta$, for every Carath\'eodory solution $t\mapsto y(t)$, defined for $t\in [0,T]$, of the perturbed Cauchy problem
\begin{equation}\label{eq:system_perturb}
\dot y=f(y,U(y+\zeta))+ d\,,\qquad\qquad y(0)=x_o\,,
\end{equation}
one has
$$
y(t)\in S\qquad\forall~t\in [0,T]\,,\qquad\qquad\mbox{ and }\qquad\qquad
y(T)\in\Sigma^\delta\,.
$$
\end{theorem}

We do not include the explicit proof of this theorem, since it follows very closely the one of Theorem~3.4 in~\cite{AB2}, by combining Theorem~\ref{thm:practical_stab_constr} with the general robustness results for patchy vector fields in presence of impulsive perturbations, proved by Ancona and Bressan in~\cite{AB2}.

\subsection{Conclusions and open problems}

In this paper we have positively solved the problem of practical stabilization of a constrained dynamics through feedback controls. As expected, the control is in general discontinuous but it is possible to select the discontinuity in a suitable way so to obtain a patchy feedback control, which in turn ensures the existence of Carath\'eodory solutions of the closed loop system for positive times, and the robustness of the feedback with respect to both inner and outer disturbances.

The problem that remains open and that is currently under investigation is the existence of nearly optimal patchy feedbacks for a constrained dynamics. In the unconstrained case, it is known that nearly optimal patchy controls exists (see~\cite{AB5, BP2}), but in the constrained case the only available result is the one contained in~\cite{CRS}, involving general discontinuous controls and Euler solutions. In view of the further robustness properties enjoyed by patchy controls, it would be of interest to provide a similar construction in terms of patchy feedbacks and Carath\'eodory trajectories.

\bigskip

\n {\bf Acknowledgement.} This work has been supported by the European Union Seventh Framework Programme [FP7-PEOPLE-2010-ITN] under grant agreement n.264735-SADCO. The author warmly thanks Prof. Fabio Ancona for his careful reading of the manuscript, his precise comments and for the repeated discussions on  the various aspects of the construction.

\end{document}